\documentclass[a4paper,12pt]{article}
\usepackage[dvips]{graphicx}

\font\Sets=msbm10
\def\Natural{\hbox{\Sets N}}

\def\bb{\begin{equation}}
\def\ee{\end{equation}}
\def\ban{\begin{eqnarray}}
\def\ean{\end{eqnarray}}
\def\ba{\begin{eqnarray*}}
\def\ea{\end{eqnarray*}}
\def\const{\hbox{const}}
\def\res{\hbox{res}}
\def\pt{\partial}

\def\a{\alpha}
\def\b{\beta}
\def\D{\Delta}
\def\d{\delta}
\def\G{\Gamma}
\def\g{\gamma}
\def\ve{\varepsilon}
\def\k{\kappa}
\def\l{\lambda}
\def\m{\mu}
\def\n{\nu}

\def\O{\Omega}

\def\s{\sigma}
\def\S{\Sigma}
\def\t{\tau}
\def\th{\theta}
\def\z{\zeta}
\newtheorem{demo}{Remark}

\title{\bf Asymptotic description of nonlinear
resonance \footnote{This work was supported by grants RFBR
(00-01-00663, 00-15-96038) and INTAS (99-1068)}}

\author{O.M. Kiselev\footnote{Institute of Mathematics
of Ufa  Sci. Centre of  RAS;\ e-mail:ok@ufanet.ru} \and S.G.
Glebov\footnote{Ufa State Petroleum Technical University;\
e-mail:gloomy69@ufanet.ru}}
\date{May 2, 2001}

\begin{document}
\maketitle

\begin{abstract}
We study a hard regime of stimulation  of  two-frequency
oscillations in the main resonance equation with a fast
oscillating external force: $ \ve i \psi' + |\psi|^2\psi =
\exp\big(it^2/ (2\ve)\big),\,\, 0<\ve\ll1$. This phenomenon is
caused by resonance between an eigenmode and the external
force. The asymptotic solution before, inside and after the
resonance layer  is studied in detail and matched.
\end{abstract}

\section{Introduction}
\par
In this paper we investigate the hard mode of stimulating of
two phase oscillations. We study this phenomenon  in
asymptotic solution of ordinary nonlinear differential
equation under fast oscillating external force:
\bb
\ve i \psi' + |\psi|^2\psi = \exp\bigg({{it^2}\over
{2\ve}}\bigg), \label{sch}
\ee
where $0<\ve\ll1$ -- a small parameter.
\par
The solution of the equation (\ref{sch}) constructed in this
paper oscillates with the frequency of the perturbation when
$t>t_*$. An amplitude of  oscillations changes slowly. This
system changing corresponds to the stimulated oscillations of
nonlinear equations. The frequency of the oscillations changes
with time. From the other hand, the eigen-frequency of the
nonlinear equation depends on the amplitude. At certain moment
$t=t_*$ the eigen-frequency of the oscillations becomes equal
to the  external force frequency. A resonance takes place in
the system. This leads to hard loss of stability of the
eigenmode of the oscillations and the system gets out of the
resonance. The asymptotic solution of equation (\ref{sch})
becomes two-phase. One phase relies to the  oscillations
stimulated by the perturbation force and another one
corresponds to the oscillations
 stimulated by the transition over the resonance.
\par
To explain the nature of the solution bifurcation and to construct
the asymptotic solution it is  more convenient to investigate
the equation for the  amplitude of the oscillations
$U=\psi\exp\big(-it^2/(2\ve)\big)$
\begin{equation}
\ve i  U' + |U|^2 U - tU = 1. \label{sh}
\end{equation}
\par
If we begin to study the amplitude $U$ and consider the equation
(\ref{sh}) then constructing of the asymptotic solution is equivalent
 to investigation of the bifurcation of the slowly varied equilibrium of
the equation (\ref{sh}). The theory of the bifurcation usually
investigates the behavior of a solution depending on external
parameter which is not connected with a variable of a
differential equation \cite{Andr,Arn}.
\par
The bifurcation of the equilibrium with varying  parameter is
well known in physics \cite{Z-S}. This bifurcation is usually
illustrated by loss of stability of the  equilibrium of
mechanical system, which is described by the  second order
differential equation with slowly varying coefficients. In
this case the hard loss of stability takes place. The numeric
evaluations give the picture:
\begin{figure}[h]
\includegraphics[angle=0, bb=0 0 300 200]{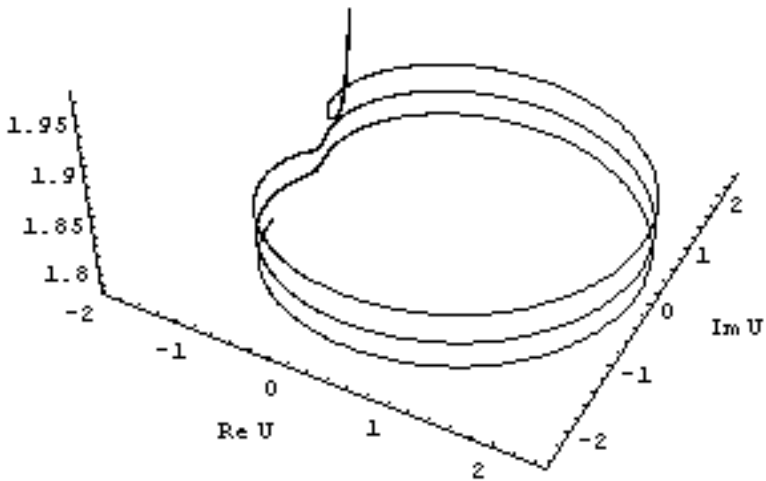}
\caption{The bifurcation solution of the equation (\ref{sh}).}
\end{figure}
\par
A description of the internal asymptotic structure of the
solution in the case of the hard loss of the stability is
connected with a construction of the asymptotics for the
solution of the Painlev\'e-1 equation \cite{Hab1}. Besides the
solution has a complicated asymptotic structure in the
transition layer where the main term of the asymptotics is
defined by four various expansions of different types. There
are: a special solution of the Painlev\'e-1 equation, a
sequence of separatrix solutions of the homogeneous equation
with "frozen" coefficients, a sequence of solutions of the
Weierstrass equation with the parameter $g_2=0$ and a sequence
of solutions of the Weierstrass  equation with the parameter
$g_2\not=0$.
\par
Two of the first internal layers were found and studied in
\cite{Hab1}. Later, in \cite{ok1}, \cite{Dim-Hab} the
asymptotic structure of the transition layer was studied in
detail. It was shown that actually it is necessary to study
the sequence of the separatrix solutions of the homogeneous
equation with "frozen" coefficients. The third internal layer
was found in \cite{Dim-Hab}.
\par
However, these three layers  are not enough for a passage
through the bifurcation interval. There exists a layer which
is defined by a solution of the  Weierstrass  equation with
the parameter $g_2\not=0$.  This layer is found and studied
in this work due to the matching of all asymptotic solutions
before inside and after the bifurcation layer. The another
new result of this paper is the calculating of internal
variables in a sequence of layers, which are connected with
the separatrix solutions of the homogeneous equation with
"frozen" coefficients. For the Painlev\'e-2 equation with a
small parameter at a derivative these calculations have been
made in the work of one of the authors of this paper  in
\cite{ok2}.
\par
Among other works devoted to the investigation of the
asymptotics solutions for nonlinear differential equations
with variable coefficients it is necessary to note the
investigations of passage through a separatrix in \cite
{Timof} -\cite {Hab-Ho}, where there were no confluent of the
slowly varying equilibrium. In particular, in work
\cite{Hab-Ho} the passage through a separatrix of equation
(\ref{sh}) with constant coefficients and small dissipation
was studied. In the case considered in this work the solution
of the equation for the amplitude also passes through a
separatrix, but at more complicated confluent saddle-center
equilibrium.
\par
Here, the asymptotic solution for the equation (\ref{sch}) in
the domain where the amplitude varies slowly ($t>t_*$) is
constructed by perturbation theory method and in domain where
the amplitude fast oscillates --- by the Krylov-Bogolyubov
method \cite{K-B, Kuz}. Thus, we  reproduce the elegant
formulas for the asymptotic  solutions of the equation
(\ref{sch}) obtained in the work  \cite{Bour-Hab2}. All
asymptotics are matched \cite{I}.
\par
The structure of this paper is following. In the second
section the problem is formulated. In the third one, the
obtained results are represented. The fourth section is
devoted to a construction of an algebraic asymptotic solution
for the equation (\ref{sh}). The fifth section includes an
investigation of a bifurcation layer and studying of the
sequences of the internal expansions. In the sixth section
the fast oscillating solution of the equation (\ref{sh}) is
constructed and its asymptotics at the approach to a
bifurcation point is written.

\section{Studied problem and its justification}
\subsection{Justification}
\par
The general solution of the equation (\ref{sch}) generally
speaking is unknown. But we can construct a set of formal
asymptotic solutions on a small parameter $\ve$. The simplest
kind of the solutions is the solutions oscillating with the
frequency of the external force. To obtain this kind of
solutions it is necessary to proceed the amplitude equation
(\ref{sh}). Further following the natural suggestion about
boundedness of the derivative in the equation (\ref{sh}) we
can obtain a nonlinear algebraic equation for the main term of
asymptotics $\stackrel{0}{U}(t)$:
\bb
|\stackrel{0}{U}|^2\stackrel{0}{U} - t\stackrel{0}{U} = 1.
\label{alg}
\ee
The number of the roots of this algebraic equation depends on
a parameter $t$. The roots can be written explicitly. There
exist a value of the parameter $t$ equals to
$t_*=3(1/2)^{2/3}$ so that the equation (\ref{alg}) has three
real roots at $t>t_*$. At $t=t_*=3(1/2)^{2/3}$ there is one
simple root and one multiple root  $U_*=-(1/2)^{1/3}$. At
$t<t_*$ the equation (\ref{alg}) has the alone root.
\par
Let us consider the domain $t>t_*$ where there are three
different roots of the equation (\ref{alg}).  Denote them
$U_k$, and $U_3(t)<U_2(t)<U_1(t)$. These roots correspond to
slowly varying equilibriums of the equation (\ref{sh}). Two
of the equilibriums are stable centers. The formal asymptotic
solutions with the main terms $U_1(t),\, U_2(t)$ correspond
to these stable equilibriums. The third equilibrium is a
saddle. There exists a parameter value $t=t_*$ at which one
of the centers coalesces with the saddle. At that moment the
saddle-center bifurcation takes place. At $t<t_*$ there
exists only one slowly varying equilibrium.
\par
This bifurcation may by explained on the example of an
autonomous  equation with a "frozen" coefficient $T$:
$$
iV'+(|V|^2-T)V=1.
$$
\begin{figure}[h]
\includegraphics[angle=0, bb=0 0 380 100]{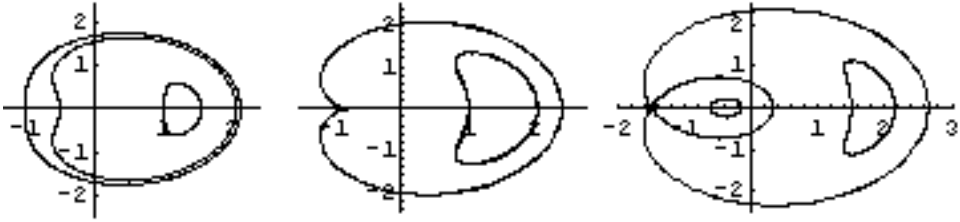}
\caption{The phase portraits for the equations with the
"frozen" coefficient  at $T<t_*$, $T=t_*$ and $T>t_*$.}
\end{figure}
On the figure 2 one can see the phase plane of this equation.
On the left picture one equilibrium exists. On the middle
picture one can find two equilibrium positions and at last on
the right picture one can find three equilibrium positions.
\par
Choosing $U_k$ as the main term of the asymptotics in the
domain $t>t_*$ one can construct three different formal
asymptotic solutions in the form:
\bb
U(t;\ve) = \sum_{n=0}^{\infty}\ve^n\stackrel{n}{U}(t). \label{alas}
\ee
Each of them corresponds to the one of the slowly varying
equilibrium. The formal solution with the main term $U_1(t)$
doesn't change the structure when the variable $t$ changes and
we do not investigate this case. Solution based on $U_3(t)$ is
unstable on the whole interval $t>t_*$ and it is also out of
our consideration.  We will investigate the asymptotic
solution based on $U_2(t)$ and the behavior of this asymptotic
solution after the saddle-center bifurcation.

\subsection{Statement of the problem}

In this work we solve the problem on constructing of the
formal asymptotic solution of the equation (\ref{sch})  in the
interval $t\in[t_*-C,t_*+C]$ where $ C=\const>0$ uniform on
$\ve$ . We suppose that the solution in the domain $t>t_*$ has
the form $$ \psi(t,\ve)=\exp\bigg({it^2\over2\ve}\bigg)
\sum_{n=0}^{\infty}\ve^n\stackrel{n}{U}(t),\quad \ve \to 0, $$
where $\stackrel{0}{U}(t)=U_2(t)$.

\section{Result}
\par
The combined asymptotic solution of the equation (\ref{sch})
in an interval  $t\in[t_*-C,t_*+C]$ is constructed. The
solution has the different asymptotic structure in different
parts of the interval. The main result of this work is
following:  firstly, the domain where constructed combined
solution is valid completely covers the interval
$[t_*-C,t_*+C]$ and secondly, all asymptotics are matched. In
this section we only represent the form of asymptotics and
leader terms. The explicit formulas can be found in the
corresponding sections of the paper.
\par
In the domain $(t-t_*)\ve^{-4/5} \gg 1$ the asymptotic solution
has the form: $$ \psi(t,\ve)= \exp\bigg({it^2\over2\ve}\bigg)
\sum_{n=0}^\infty\ve^{n}\stackrel{n}{U}(t),\quad \ve \to 0. $$
The leader term of the asymptotics is equal to the middle root of
algebraic equation (\ref{alg}):  $\stackrel{0}{U}(t)=U_2(t)$. The
correction terms $\stackrel{n}{U}(t)$ are algebraic functions of
$t$.
\par
In the domain $|t-t_*|\ll 1$ the asymptotic solution is defined
by four various expansions of different types. First of them is:
\bb \psi=\bigg(U_* +
\ve^{2/5}\sum_{n=0}^\infty\ve^{2n/5}\stackrel{n}{\alpha}(\tau) +
i\ve^{3/5}\sum_{n=0}^\infty\ve^{2n/5}\stackrel{n}{\beta}(\tau)
\bigg)\exp\bigg({it^2\over2\ve}\bigg),\quad \ve \to 0, \label{as2}
\ee where $U_*$ is multiple root of the cubic equation
(\ref{alg}) at $t=t_*$. The coefficients of this expansion
depend on new scaling time $\t=(t-t_*)\ve^{-4/5}$. The leader
term of the asymptotics and the corrections are defined by
their asymptotics as $\t\to-\infty$ uniquely. In particular,
the leader term of asymptotics $\stackrel{0}{\a}(\t)$ is a
special solution of the Painlev\'e-1 equation: $$
\stackrel{0}{\alpha}{}^{\prime\prime} -
3\stackrel{0}{\alpha}\!\!^2+\tau = 0. $$ with the given
asymptotics as $\t\to-\infty$: $$ \stackrel{0}{\alpha}(\tau) =
\sum_{n\ge0} \alpha_n \tau^{-{{(5n-1)}\over
2}},\quad{\hbox{where}}\,\,\,\,
\a_0={1\over\sqrt{3}},\,\,\a_1={1\over24}. $$ In the domain
$\t>-\infty$ this solution has poles on the real axis of $\t$.
Denote the least of them by $\t_0$. The asymptotics
(\ref{as2}) is valid as $(\tau-\tau_0)\ve^{-1/5}\gg 1$.
\par
In the neighborhood of the bifurcation point $(\t=\t_0)$ the
coefficients of the asymptotic expansion depend on one more
fast time scale $\theta=(\t-\t_0)\ve^{-1/5}$. Denote by $$
\th_0=\th+\sum_{n=1}^\infty\ve^{n/5}\stackrel{n}{\th}_0, $$
where $\stackrel{n}{\th}_0=\const$. Then in the domain
$-\ve^{-1/5}\ll\th_0\ll\ve^{-1/10}$ the formal asymptotic
solution has the form: $$ \psi(t,\ve)
=\bigg(U_*+\stackrel{0}{w}(\th_0)+\ve^{4/5}
\sum_{n=1}^\infty\ve^{(n-1)/5}\stackrel{n}{w}(\th_0)\bigg)
\exp\bigg({it^2\over2\ve}\bigg), \quad \ve \to 0. $$ The main
term of asymptotics $\stackrel{0}{w}(\th_0)$ is the separatrix
solution of the autonomous equation: \bb
i\stackrel{0}{w}{}^\prime +
U_*\left(2|\stackrel{0}{w}|^2+\stackrel{0}{w}\!{}^2\right) +
U_*^2 \left( \stackrel{0}{w}\!^* - \stackrel{0}{w} \right) +
|\stackrel{0}{w}|^2 \stackrel{0}{w}= 0, \label{avt} \ee
namely: $$ \stackrel{0}{w}(\th_0)
={{-2}\over{(\th_0-iU_*)^2}}. $$
\par
In the domain $-\th_0\gg1$ the asymptotic solution is defined
by a sequence of two alternating asymptotics. Let us call
them by "intermediate" and "separatrix" asymptotics. To obtain
the intermediate asymptotics  let us introduce one more slow
variable:
$$ T_k=\th_{k-1}\ve^{1/6}, k=1,2,\dots. $$ An asymptotic solution
in the intermediate domain for not too large values $ k $: $$
k\ll\ve^{-1/7} $$ has the form:
$$
\psi(t,\ve)=\bigg(U_*+\ve^{1/3}\sum_{n=0}^{\infty}\ve^{1/30}
\stackrel{n}{A}_k(T_k)
\,+\,
$$
$$
i\ve^{1/2}\sum_{n=0}^{\infty}\ve^{1/30}
\stackrel{n}{B}_k(T_k)\bigg)\exp\bigg({it^2\over2\ve}\bigg),\quad
\ve \to 0. $$ The leader term of the asymptotics satisfies to
the equation: $$
\stackrel{0}{A}_k\!\!\!''+3\stackrel{0}{A}_k\!\!\!^2=0 $$ and
can be expressed by the Weierstrass function with the
parameter $g_2=0$:
\bb
\stackrel{0}{A}_k=-2\wp(T_k;0,g_3(k)),\,\,\,\hbox{where}
\,\,\, g_3(k)={1\over56}\big(g_3(k-1)+\pi/2\big),\,\,\,
g_3(0)={a_4\over56}. \label{A-0-0}
\ee
The constant $a_4$ is
the coefficient as $(\t-\t_0)^4$ in the Laurent expansion of
$\stackrel{0}{\a}(\t)$.
\par
The Weierstrass $\wp$-function has a real period $\O_k$ and
has poles on the real axis at $T_k=0$ and $T_k=\O_k$. The
intermediate expansion with the leader term (\ref{A-0-0}) is
valid in the domain between the poles as
$$
-\ve^{-1/6}T_k\gg1,\quad \ve^{-2/15}(T_k+\O_k)\gg1.
$$
\par
At the large values of $k$ the intermediate asymptotics are
constructed in the form
$$
\psi(t,\ve)=\bigg(U_*+\ve^{1/3}\sum_{n=0}^{\infty}\ve^{n/6}
\stackrel{5n}{A}_k(T_k)
\,+\,
$$
$$
i\ve^{1/2}\sum_{n=0}^{\infty}\ve^{n/6}
\stackrel{5n}{B}_k(T_k)\bigg)\exp\bigg({it^2\over2\ve}\bigg),\quad
\ve \to 0.
$$
The leader term of asymptotics satisfies to the
equation: $$
\stackrel{0}{A}_k''+3\stackrel{0}{A}_k\!\!\!^2=\l_k.
$$ $$ \l_k(\ve)=\ve^{1/6}\bigg(\sum_{j=1}^{k}\O_j
\,+\,\sum_{n=1}^\infty\ve^{(n-1)/30}
\sum_{j=1}^k\stackrel{n}{x}_j\!\!\!^+\bigg). $$ The main term
of the asyptotics is: \bb
\stackrel{0}{A}_k(T_k)=-2\wp(T_k,\l_k/2,g_3(k,\ve)),
\label{A-0-lk} \ee where $$
g_3(k,\ve)=\stackrel{0}{g}_3(k)+\sum_{n=1}^\infty
\ve^{n/30}\stackrel{n}{g}_3(k). $$ The intermediate expansion
with the leader term (\ref{A-0-lk}) is valid in the domain
between the poles of the Weierstrass function as $$
-\ve^{-1/6}T_k\gg1,\quad \ve^{-2/15}(T_k+\O_k)\gg1. $$
\par
The separatrix expansions are valid in  a small neighborhood
of the  Weierstrass function poles. Denote: $$
\th_k=\big(T_k+\O_k-{1\over4}\sum_{n=1}^\infty\ve^{n/30}
\stackrel{n}{x}_k\!\!\!^+\big)\ve^{-1/6},\quad k=1,2,\dots. $$
When $$ |\th_k|\ve^{1/6}\ll1 $$ the formal asymptotic
solution of equation (\ref{sch}) has the form: $$
\psi=\bigg(U_*+\stackrel{0}{W}(\th_k)+\ve^{4/5}\sum_{n=1}^\infty
\ve^{(n-1)/30}\stackrel{n}{W}(\th_k)\bigg)\exp\bigg({it^2\over2\ve}\bigg)\quad
\ve \to 0. $$ The leader  term of the asymptotics
$\stackrel{0}{W}(\th_k)$ is a separatrix solution of the
autonomous equation  (\ref{avt}): $$ \stackrel{0}{W}(\th_k)
={{-2}\over{(\th_k-iU_*)^2}}. $$
\par
The  sequence of the alternating intermediate expansions and
separatrix asymptotic ones is valid as
$\ve^{-1/6}(t_*-t)\ll1$.
\par
In the domain $(t_*-t)\ve^{-2/3}\gg1$ the asymptotic solution
becomes two-phase. The amplitude of the stimulated
oscillations in the solution of (\ref{sch}) oscillates fast.
The form of the solution is:
$$
\psi=\bigg(\stackrel{0}{U}(t_1,t,\ve) + \ve
\stackrel{1}{U}(t_1,t,\ve) + \ve^2 \stackrel{2}{U}(t_1,t,\ve)
\bigg)\exp\bigg({it^2\over\ve}\bigg),
$$
where $t_1$ is a new fast variable $t_1=S(t)/\ve+\phi(t)$.
The main term of the asymptotics $\stackrel{0}{U}$ lies on the
curve $\G(t)$:
$$
{1\over 2}|y|^4 -t |y|^2 - (y + \bar y) = E(t),
$$
and satisfies to the Cauchy problem for the equation
$$
iS^{\prime}\pt_{t_1} \stackrel{0}{U} = \pm
\sqrt{2\stackrel{0}{U}\!^3 + (2E(t)+t^2) \stackrel{0}{U}\!^2 +
2t\stackrel{0}{U}+1},
$$
with an initial condition  $\stackrel{0}{U}|_{t_1=0}=u_0$, such,
that $\Im(u_0)=0,\quad \Re(u_0)=\min_{y\in
\G(t)}\big(\Re(y)\big)$. The function $S(t)$ is a solution
for the Cauchy problem
$$
iS'\int_{\G(t)}{dy\over\sqrt{3y^3+(2E+t^2)y^2+2ty+1}}=T, S|_{t=0}=0,
$$
where  $T=\const>0$. The function $E(t)$ is the solution of
the transcendental equation
$$
i\int_{\G(t)}u^*du=\pi.
$$
The phase shift $\phi$ is defined by two constant. One of
them is an "initial" condition $\phi(t_*)=\phi_0$ and another one
is a constant $\phi_1$ in the equation in for  $\phi$:
$$
{\phi'\over \pt_E S}\pt_E I =\phi_1.
$$
\par
\begin{demo} In our work  constants $\phi_0$ and $\phi_1$ are not
calculated. It means that the phase shift of the fast
oscillating asymptotics remains to be undefined. Its
evaluation within the framework of our approach requires an
evaluation on an explicit form of the asymptotics of the
higher corrections in the  internal expansions at their
matching with the fast oscillating asymptotics.
\end{demo}

\section{Construction of an algebraic asymptotics}

In this section the formal asymptotic solution for the
equation (\ref{sh}) has been constructed  in the form of
series on integer powers of the small parameter $\ve$. We show
 that this expansion is valid in the
interval $(t - t_*)\ve^{-4/5} \gg 1$.
\par
The coefficients of this asymptotics are defined from the
recurrent sequence of algebraic equations. These equations
can be obtained by substituting of the series (\ref{alas}) in
the equation (\ref{sh}) and by collecting of the terms under
the same powers of the small parameter. In particular, the
algebraic equation (\ref{alg}) for the leader term of the
asymptotics can be obtained from the relation under $\ve^0$.
\par
Further  we will use the root $U_2(t)$ as the main term
$\stackrel{0}{U}(t)$ of the asymptotics (\ref{alas}). It means
that we investigate the formal asymptotics (\ref{alas})
corresponding to the stable center which converges with the
saddle at $t_*$.
\par
The relation under $\ve^1$ give us the equations for real and
imaginary parts of the first correction term
$\stackrel{1}{U}=u_1 +i v_1$
\bb
[3\stackrel{0}{U}\!^2-t]u_1=0,\quad
[\stackrel{0}{U}\!^2-t]v_1=-\pt_t \stackrel{0}{U}. \label{1st}
\ee
The expression $[3\stackrel{0}{U}\!^2-t]$  vanishes at
$t=t_*$. So it is easy to obtain that the even corrections
are real, and the odd corrections are imaginary.
\par
The explicit form of the first correction is:
\bb
\stackrel{1}{U} = - {i\,{\stackrel{0}{U}}\over
{[3\stackrel{0}{U}\!^2-t][\stackrel{0}{U}\!^2-t]}}.
\ee
This expression has a singularity at $t=t_*$. The order of
the singularity increases with respect to the number of the
corrections of the asymptotics as $t \to t_*$. In the
neighborhood of the point $t_*$ the coefficients of
asymptotics (\ref{alas}) have the form:
\par
in the case of $n=0,1$:
\bb
\stackrel{n}{U}(t)=(t-t_*)^{-n/2}\sum_{k=0}^\infty\stackrel{n}{U}_k
(t-t_*)^{k/2}, \label{expU01} \ee
where $\stackrel{0}{U}_0=U_*$,
$\stackrel{0}{U}_1=1/\sqrt{3}$;
\par
in the case of $m=2n,\,\, m=2n+1$, $n\in\Natural$:
\bb
\stackrel{m}{U}(t)=(t-t_*)^{(1-5m)/2}\sum_{k=0}^\infty\stackrel{m}{U}_k
(t-t_*)^{k/2}. \label{expUnm}
\ee
\par
The increasing of the singularity of the corrections in
(\ref{alas}) makes this expansion unusable as $t\to t_*$. The
singularities of the coefficients define the domain  of the
suitability for the formal asymptotic solution (\ref{alas})
as $(t - t_*)\ve^{-4/5} \gg 1$.

\section{Expansions in bifurcation layer}

In this section the formal asymptotic solution for equation
(\ref{sh}) has been constructed  in the domain $|t-t_*| \ll
1$. Usually the  asymptotics of this type is called internal
ones \cite{I}.

\subsection{Initial interval}
\subsubsection{The first internal expansion --- Painlev\'e-1 interval}

In the neighborhood of the singularity we will use a new
scaling variable $\tau=(t-t_*)\ve^{-4/5}$. This new variable
$\tau$ is defined by the structure  of external asymptotic
solution (\ref{alas}) singularities.
\par
We will construct the solution of the original equation (\ref{sh})
in the form:
\bb
U=U_* + \ve^{2/5}\alpha(\tau,\ve) + i\ve^{3/5}\beta(\tau,\ve),
\label{infs}
\ee
\par
The equations for $\a$ and $\b$ are
\begin{eqnarray}
\a'+(U_*^2-t_*)\b=-\ve^{2/5}2U_*\a\b\,-\,\ve^{4/5}(\a^2-\t)\b\,-\,\ve^{6/5}\b^3,
\nonumber\\
\b'-3U_*\a^2+U_*\t=\ve^{2/5}(-\a^3-U_*\b^2+\a\t)\,+\,\ve^{4/5}\a\b^2.
\label{eq20}
\end{eqnarray}
\par
The asymptotics of the functions $\a$ and $\b$ as
$\t\to\infty$ is known. It can be obtained by substituting
(\ref{expU01}) and (\ref{expUnm}) in formula (\ref{alas}) and
by decomposing of this expression in the terms of the scaling
variable $\t$:
\begin{eqnarray}
\alpha(\tau,\ve) = {1\over\sqrt{3}}\sqrt{\tau}
+\sum_{n=1}^{\infty} \tau^{(1-5n)/2}\stackrel{2n}{U}_0 +
\ve^{2/5} \sum_{n=0}^{\infty}\stackrel{2n}{U}_1
\tau^{(2-5n)/2} +\nonumber\\
\ve^{4/5}\sum_{n=0}^{\infty}\stackrel{2n}{U}_2
\tau^{(3-5n)/2}+\dots.
\label{asst}
\end{eqnarray}
\begin{eqnarray}
\b(\t,\ve)={1\over\sqrt{\t}}\stackrel{1}{U}_0+
\sum_{n=1}^{\infty}\stackrel{2n+1}{U}\!\!\!\!_0\,\,\t^{-(1+5n)/2}\,+\,
\ve^{2/5}\sum_{n=0}^\infty
\stackrel{2n+1}{U}\!\!\!\!_1\,\,\tau^{-(2-5n)/2} +\dots. \label{asbt}
\end{eqnarray}
\par
It is convenient  to construct the solution of the equations
(\ref{eq20}) in the form of formal series on powers of the
small parameter $\ve$:
\bb \alpha(\tau,\ve)=
\sum_{n=0}^{\infty}\ve^{2n/5}\stackrel{n}{\alpha}(\tau),\quad
\beta(\tau,\ve)=
\sum_{n=0}^{\infty}\ve^{2n/5}\stackrel{n}{\beta}(\tau). \label{as20}
\ee
The leader correction terms are defined by the system of ordinary
differential equations
\begin{eqnarray}
\stackrel{0}{\alpha}\!' + (U_*^2 - t_*)\stackrel{0}{\beta}=0,\\
\stackrel{0}{\beta}\!'- 3U_*\stackrel{0}{\alpha}\!\!^2+\tau U_* =0.
\end{eqnarray}
Let us differentiate the first equation and  take into account
that $U_*(U_*^2-t_*)=1$ then we obtain the Painlev\'e-1
equation for $\stackrel{0}{\alpha}$
\bb
\stackrel{0}{\alpha}{}^{\prime\prime} -
3\stackrel{0}{\alpha}\!\!^2+\tau = 0. \label{p1}
\ee
Matching the asymptotics (\ref{infs}), (\ref{as20}) and
asymptotics (\ref{alas}) as $t\to t_*+0$ we obtain the
asymptotic structure of $\stackrel{0}{\alpha}(\tau)$ as $\tau
\to +\infty$
\bb
\stackrel{0}{\alpha}(\tau) = \sum_{n\ge0} \alpha_n
\tau^{-{{(5n-1)}\over 2}},\quad{\hbox{where}}\,\,\,\,
\a_0={1\over\sqrt{3}},\,\,\a_1={1\over24}. \label{asp1}
\ee
The solution of the Painleve-1 equation with the asymptotics
of (\ref{asp1}) was investigated in \cite{H-S,G-L}.
\par
The coefficients of the formal expansion (\ref{as20}) for the function
$\b$ can be obtained from the first equation of the system (\ref{eq20}).
In particular:
$$
\stackrel{0}{\b}={\stackrel{0}{\a}\!'\over2U_*^2},\quad
\stackrel{1}{\b}={1\over2U_*^2}\bigg(\stackrel{1}{\a}\!'+2U_*\stackrel{0}{\a}
\stackrel{0}{\b}\bigg),
$$
$$
\stackrel{2}{\b}={1\over2U_*^2}\bigg(\stackrel{2}{\a}\!'+
\stackrel{0}{\a}\!^2\stackrel{0}{\b}+2U_*\stackrel{1}{\a}\stackrel{0}{\b}+
2U_*\stackrel{0}{\a}\stackrel{1}{\b}-\t\stackrel{0}{\b} \bigg).
$$
\par
The equation for the first correction term $\stackrel{1}{\alpha}(\tau)$
has the form
$$
\stackrel{1}{\alpha}{}^{\prime\prime} +
6U_*(U_*^2-t_*)\stackrel{0}{\alpha} \stackrel{1}{\alpha}
+2U_*(U_*^2-t_*)\stackrel{0}{\beta}\!{}^2 +
(U_*^2-t_*)\stackrel{0}{\alpha}\!{}^3
-(U_*^2-t_*)\stackrel{0}{\alpha}\tau =
$$
\bb
-2U_*\stackrel{0}{\alpha}{}^\prime\stackrel{0}{\beta} -
 2U_*\stackrel{0}{\alpha}\stackrel{0}{\beta}{}^\prime. \label{al1}
\ee
The asymptotics of a solution for this equation is uniquely
defined by matching of the asymptotics of the external
expansion (\ref{alas}) and internal expansion (\ref{infs}),
(\ref{as20}):
\bb
\stackrel{1}{\alpha} = \stackrel{0}{U}_2\tau -
\stackrel{2}{U}_6\tau^{-3/2} + O(\tau^{-4}), \quad \tau \to +\infty
\label{a1-as+}
\ee
A simple calculation gives us the domain of a suitability of
that asymptotic expansion. The representation (\ref{as20})  of
the solution is valid as $\tau \ll \ve^{-4/5}$.  This result
is similar to \cite{Hab1}.
\par
It is known that the solution of the Painlev\'e-1 equation
with the asymptotics (\ref{asp1}) has poles at points $\tau_k$
\cite{Kit}. Let us denote the first pole by $\tau_0$ and call
it the point of the bifurcation.
\par
The solution of the Painlev\'e-1 equation is represented by
converging series in the neighborhood of the point $\tau_0$
\cite{G-L}:
 $$
\stackrel{0}{\alpha}(\tau)=-{2\over{(\tau-\tau_0)^2}} -
{\tau_0\over10}(\tau-\tau_0)^2 + a_4(\tau-\tau_0)^4 +
\sum_{k=6}^\infty\stackrel{0}{\a}_k(\tau-\tau_0)^k, $$
the constants $\tau_0$ and  $a_4$ are parameters of the solution.
\par
The first correction term  $\stackrel{1}{\a}$ is defined by
the equation  (\ref{al1}) and asymptotics (\ref{a1-as+}). In the neighborhood
of the bifurcation point $\t_0$ as $\t>\t_0$ the function
$\stackrel{1}{\a}$ is represented in the form:
\bb
\stackrel{1}{\a}=\stackrel{1}{\a}_c(\t)+\stackrel{1}{a}\!^+
\a_1(\t)+\stackrel{1}{b}\!^+\a_2(\t). \label{a1-as-bif}
\ee
Where  $\stackrel{1}{\a}_c(\t)$ is the partial solution
of the nonhomogeneous equation for the function
$\stackrel{1}{\a}(\t)$ with the asymptotics:
$$
\stackrel{1}{\alpha}_c(\tau) =
\stackrel{1}{\a}_{-4}(\tau-\tau_0)^{-4} +
\sum_{k=0}^\infty\stackrel{1}{\a}_k(\t-\t_0)^k, \quad
\stackrel{1}{\a}_{-4}={1\over {U_*^7}},\quad \stackrel{1}{\a}_4=0.
$$
It is convenient to use the partial solution with the
coefficients as $(\t-\t_0)^{-3}$ and $(\t-\t_0)^4$ equal to zero.
\par
The constants $\stackrel{1}{a}\!^+$ and $\stackrel{1}{b}\!^+$ in
the formula (\ref{a1-as-bif}) can be uniquely defined by
asymptotics of $\stackrel{1}{\a}$ as $\t\to\infty$. The functions
$\a_1(\t)$ and $\a_2(\t)$ are linear independent solutions of the
homogeneous equation for $\stackrel{1}{\a}(\t)$. These functions
are uniquely defined by their asymptotics as $\t \to \t_0$:
$$
\a_1(\t)=(\t-\t_0)^{-3}+a_1(\t-\t_0)+a_2(\t-\t_0)^2+
a_3(\t-\t_0)^3+O((\t-\t_0)^5),
$$
where $a_n=\const$;
$$
\a_2(\t)=(\t-\t_0)^4+O((\t-\t_0)^8).
$$
\par
The next correction terms are represented by the following
formula:
$$
\stackrel{n}{\a}(\t)=\stackrel{n}{\a}_c(\t)+\stackrel{n}{a}\!^+\a_1(\t)
+\stackrel{n}{b}\!^+\a_2(\t).
$$
The partial solution $\stackrel{n}{\a}_c(\t)$ of the equation has
the asymptotics as $\t \to \t_0$:
$$
\stackrel{n}{\a}_c(\t)=\stackrel{n}{\a}_{-2n}(\t-\t_0)^{-2n}+\dots,
$$
where dots mean the terms of the next orders on $(\t-\t_0)$.
\par
It is convenient for us to denote  the partial solution which
does not contain terms of the order of $(\t-\t_0)^{-3}$ and
$(\t-\t_0)^{4}$ in the asymptotics as $\t \to \t_0$ by
$\stackrel{n}{\a}_c(\t)$.
\par
By analysis of written asymptotics we obtain that the formal
asymptotic solution (\ref{infs}) is valid as
$(\tau - \tau_0)\ve^{-1/5} \gg 1$.

\subsubsection{The second internal expansion, the neighborhood
of the bifurcation point}

In this section we construct the formal asymptotic solution
for the equation (\ref{sh}), which is valid in the
neighborhood of the bifurcation point $|\tau - \tau_0| \ll 1$.
\par
The domain where the first asymptotic expansion is valid
determines  the second internal scale. We will use a new
scaling variable $\th=(\tau - \tau_0)\ve^{-1/5}$ in the small
neighborhood of the point $\tau_0$.
\par
We change the unknown function in the equation (\ref{sh}) and pass to
the function $w=U-U_*$.
\par
Substituting the expressions for $w$ and $\th$ in the original
problem we obtain the equation for $w(\th,\ve)$
\begin{eqnarray}
i w' +2U_*^2w+2U_*|w|^2+|w|^2w+U_*^2 w^* +U_* w^2- \nonumber\\
-\ve^{4/5}\tau_0 \big(U_*+w\big)-\ve\th\big( U_*+w\big)
=0.\label{w-main}
\end{eqnarray}
\par
The asymptotics of the function $w(\th,\ve)$ as $\th\to
+\infty$ can be obtained by expansing in series of the
function $U=U_*+\ve^{2/5}\a(\t,\ve)+i\ve^{3/5}\b(\t,\ve)$ on
variable $\th$:
\begin{eqnarray}
w(\th,\ve)=\bigg({1\over U_*^3\th^2}+{i\over2U_*^5\th^3} +\dots\bigg)
+\ve^{1/5}\stackrel{1}{a}\!^+\bigg(\th^{-3}+\dots\bigg)+
\nonumber\\
\ve^{3/5}\stackrel{2}{a}\!^+\bigg(\th^{-3}+\dots\bigg) +
\ve^{4/5}\bigg({2 U_8^3\t_0\over10}\th^{2}+\dots\bigg)+\dots.
\label{w-as-right}
\end{eqnarray}
In this formula dots mean the terms of the lesser order (in
brackets with respect to $\th$ and out of them with respect
to $\ve$).
\par
In this section the formal asymptotic solution for the
equation (\ref{w-main}) with the given asymptotics
(\ref{w-as-right}) as $\ve\to0$ has been constructed. We
indicate the domain where this formal solution for the
equation (\ref{w-main}) is valid.
\par
The first step of the formal solution construction for the
equation (\ref{w-main}) is a choice  of an asymptotic
sequence. This choice for the external and the first internal
expansions was defined by the equation. In this case the
asymptotics as $\th\to+\infty$ gives us the asymptotic
sequence on $\ve^{1/5}$. But the equation shows, that the
first order correction should satisfy to the homogeneous
equation. In this case we can satisfy to the asymptotic
condition as $\th\to+\infty$ by constructing of the formal
solution in the form:
\bb
w(\th,\ve)=\stackrel{0}{w}(\th_0)\,+\,
\ve^{4/5}\sum_{n=1}^\infty\ve^{(n-1)/5}\stackrel{n}{w}(\th_0),
\label{exp-w} \quad \ve\to0
\ee
where
$$
\th_0=\th+\ve^{1/5}\stackrel{1}{\th_0}+
\ve^{3/5}\stackrel{2}{\th_0}+\dots \quad \ve\to0.
$$
The correction terms $\stackrel{n}{\th}$ are defined by the
asymptotics of the formal solution as $\th\to\infty$:
$$
\stackrel{n}{\th_0}={\stackrel{n}{a}\!\!^+\over4}. $$
\par
It is convenient to pass from the equation (\ref{w-main}) to
the equation with respect to the variable $\th_0$:
\begin{eqnarray}
iw' +(2U_*^2-t_*)w+2U_*|w|^2+|w|^2w+U_*^2 w^* +U_* w^2=\nonumber\\
\bigg(\ve^{4/5}\tau_0+\ve\th_0 -\ve^{6/5}\stackrel{1}{\th}_0
-\ve^{8/5}\stackrel{2}{\th}_0-\dots\bigg)\big( U_*+w\big),
\quad\ve\to0.
\label{w-main0}
\end{eqnarray}
\par
To construct the formal solution $w(\th_0,\ve)$ we substitute
the formal series  (\ref{exp-w}) in the equation
(\ref{w-main0}) and gather the terms of the same order on
$\ve$. Then we obtain the recurrent sequence of the equations
for the leader term and the corrections of the asymptotics.
The equation for the leader term is
\bb
i\stackrel{0}{w}{}^\prime +
U_*\left(2|\stackrel{0}{w}|^2+\stackrel{0}{w}\!{}^2\right) + U_*^2
\left( \stackrel{0}{w}\!^* - \stackrel{0}{w} \right) +
|\stackrel{0}{w}|^2 \stackrel{0}{w}= 0. \label{w0}
\ee
\par
The asymptotics of $\stackrel{0}{w}$ as  $\th_0\to\infty$ can be obtained
from the asympotics of the function $w(\th_0,\ve)$
\bb
\stackrel{0}{w}(\th_0)={1\over
U_*^3\th_0^2}+{i\over2U_*^5\th_0^3}+\dots. \label{w0-as-right} \ee
\par
To construct the leader term $\stackrel{0}{w}(\th_0)$ we take
the conservation law for the equation (\ref{w0}):
$$
U_*\left[|w|^2w^* + |w|^2w \right] + U_*^2\left[{1\over2} (w^*)^2 +
{1\over2} w^2 -|w|^2 \right] + {1\over2}|w|^4 = H.
$$
The constant $H$ for our solution is defined by asymptotics
(\ref{w0-as-right}): $H=0$. Let us solve this formula for the
function $w^*$ and substitute the representation over $w$
in the equation (\ref{w0}). Then we obtain
\bb
iw^{\prime} - 2U_*w\sqrt{-U_*w}=0. \label{mainwa}
\ee
This equation can be obviously solved
\bb
\stackrel{0}{w}(\th_0) ={{-2}\over{(\th_0-iU_*)^2}}, \label{wa}
\ee
the integration constant is defined by the asymptotics
(\ref{w0-as-right}). The leader term of the asymptotics
(\ref{exp-w}) is the separatrix solution of the equation
(\ref{w0}).
\par
The next corrections $\stackrel{m}{w},\,\,m=1,2 \dots$
satisfy to the equations
\bb i\stackrel{m}{w}{}^\prime
+\left[2U_*\stackrel{0}{w}\!^* + 2U_*\stackrel{0}{w} - U_*^2 + 2
\stackrel{0}{w} \stackrel{0}{w}\!^* \right]\stackrel{m}{w} +
\left(U_* + \stackrel{0}{w}\right)^2 \stackrel{m}{w}\!^* =
\stackrel{m}{F}. \label{mwa} \ee
In particular,
$$
\stackrel{1}{F}=\tau_0 (U_* + \stackrel{0}{w}),\quad
\stackrel{2}{F}=\theta_0(U_*+\stackrel{0}{w})\quad
\stackrel{3}{F}=-\stackrel{1}{\theta_0}(U_*+\stackrel{0}{w}),
\quad \stackrel{4}{F}=0, $$ $$
\stackrel{5}{F}=-\stackrel{2}{\theta}(U_*+\stackrel{0}{w})-
2U_*|\stackrel{1}{w}|^2-U_*\stackrel{1}{w}\!\!^2. $$
The asymptotics $\stackrel{m}{w}(\th_0)$ as $\th\to\infty$
can be easy obtained from (\ref{w-as-right}).
\par
To construct the correction terms we need the solutions of the homogeneous
equation (\ref{mwa}). One of them can be obtained by differentiation
with respect to variable  $\th_0$ of the function $\stackrel{0}{w}$
\bb
w_1(\th_0)={1\over 4}{d \over {d\th_0}}\stackrel{0}{w} =
{1\over(\th_0-iU_*)^3}.\label{w1}
\ee
\par
The formula for the second independent solution is complicated
and we do not write it. The second solution grows as $|\th|\to\infty$.
 The asymptotics as $\th\to\infty$ of this solution has
the  form
\bb w_2(\th_0)=\th_0^4+
i{1 \over {2U_*^2}}\th_0^3+\dots.\label{w2}
\ee
These independent solutions are such that their
Wronskian is
$$
{\cal W}=w_1w_2^*-w_2w_1^*=-{7i\over U_*^2}.
$$
\par
By using the expressions for  $w_1$ and $w_2$ we can write
the next corrections for the asymptotic solution with
the given asymptotics as $\th_0\to+\infty$ in an
explicit form. The asymptotics  as
$\th_0\to-\infty$ of these corrections can be evaluated.
\par
It is convenient for us to represent the higher corrections of the
asymptotics (\ref{exp-w}) in the form
 \bb
\stackrel{n}{w}(\th_0)=\stackrel{n}{w}_c(\th_0)+\stackrel{n}{X}_0\!\!^-w_1(\th_0)+
\stackrel{n}{Y}_0\!\!^-w_2(\th_0), \label{1waa} \ee where,
$\stackrel{n}{w}_c(\th_0)$ is the partial solution for the n-th
correction term and this solution does not contain terms of the order
of $\th_0^4$ and $\th_0^{-3}$ in the asymptotics as $\th_0=\infty$.
The constants $\stackrel{n}{X}_0\!\!^-$ and $\stackrel{n}{Y}_0\!\!^-$
can be evaluated by
 $$
\stackrel{n}{X}_0\!\!^-=-\res_{r=\infty}\bigg[\,{1\over
r}\,\int_{-r}^{r}\big(\stackrel{n}{F}w_2^*+
\stackrel{n}{F}\!\!^*w_2\big){d\th\over{\cal W}}\bigg], $$ $$
\stackrel{n}{Y}_0\!\!^-=-\res_{r=\infty}\bigg[\,{1\over
r}\,\int_{-r}^{r}\big(\stackrel{n}{F}w_1^*+
\stackrel{n}{F}\!\!^*w_1\big){d\th\over{\cal W}}\bigg]. $$
\par
The segment of the Laurent series for the function
$\stackrel{1}{w}_c(\th_0)$ has the form
$$
\stackrel{1}{w}_c(\th_0) = {{\tau_0 U_*^3}\over 5}\th_0^2 +
i{{\tau_0 U_*}\over 5}\th_0 - {\tau_0\over{30U_*}} +
i{{3\tau_0}\over{5U_*^3}}\th_0^{-1} + O(\th_0^{-2}), \quad \th
\to\infty. $$
The coefficients as $w_1$ and $w_2$ are
$$
\stackrel{1}{X}_0\!\!^-=0,\quad \stackrel{1}{Y}_0\!\!^-=0. $$
\par
The function $\stackrel{2}{w}(\th_0)$ as $\th_0\to-\infty$
also can be represented in the form of sum of the partial solution
of the nonhomogeneous equation and solutions of the homogeneous
equation
\bb
\stackrel{2}{w}(\th_0)=\stackrel{2}{w}_c(\th_0)+\stackrel{2}{X}_0\!\!^-w_1(\th_0)+
\stackrel{2}{Y}_0\!\!^-w_2(\th_0). \label{as2ct}
\ee
The partial solution of the nonhomogeneous equation can
be represented in the form as $\th_0\to-\infty$
 $$
\stackrel{2}{w}_c(\th_0)=-{\th_0^3\over6}+i{U_*\th_0^2\over2}+{5\th_0\over12
U_*}-i{3\over4}+
{15U_*\over16\th_0}+i\big({16U_*^2\over15}+{2\over3U_*}\big)\th_0^{-2}+O(\th_0^{-4}).
$$
The Laurent series of this solution does not contain the terms of the order
of $\th_0^4$ and $\th_0^{-3}$. The coefficients as
$w_1$ and $w_2$ are
$$
\stackrel{2}{X}_0\!\!^-=0,\quad \stackrel{2}{Y}_0\!\!^-={\pi\over 2
}$$
\par
The next correction terms of the asymptotics (\ref{exp-w}) also can be
obtained by this way. The condition
$\ve\stackrel{2}{w}/(\ve^{4/5}\stackrel{1}{w})\ll1$ gives us the domain
where the expansion (\ref{exp-w}) of the asymptotic solution is valid
$$ -\ve^{-1/5}\ll\th_0\ll\ve^{-1/10}. $$
\begin{demo}
It is easy to see that the asymptotics of (\ref{exp-w}) as
$\th_0\to\infty$ and as $\th_0\to-\infty$ are different. The
difference is the term of the order of $\ve\th^4$. It means
that additional terms appear via the passage of the
neighborhood of the bifurcation point. The leader term of
these additional terms is $\stackrel{2}{Y}\!\!^-_{0}=\ve
\pi\th^4/2$. This effect leads to the appearance of the
solution $\stackrel{2}{Y}\!\!^-_0 w_2$ of the homogeneous
equation in the asymptotics of the second corrections as
$\th_0\to-\infty$.
\end{demo}

\subsubsection{Intermediate expansion, initial interval}

In this section we investigate the behaviour of the solution
after a passage of the narrow neighborhood of the bifurcation
point $\tau_0$.
\par
 In this domain it is convenient
to introduce a new scaling variable $T_1=\th_0\ve^{1/6}$ by
analysis of the second internal expansion.
\par
The solution of the equation (\ref{sh}) we seek in the form:
\bb
U=U_* + \ve^{2/6} A(T_1,\ve) + i \ve^{3/6}B(T_1,\ve),
\label{3a}
\ee
where $A(T_1,\ve)$ and $B(T_1,\ve)$ are real functions.
\par
The equation (\ref{sh}) can be  written in terms of these new
variables  in the form:
\begin{eqnarray}
A'-2U_*^2 B=-\ve^{1/3}2U_*AB-\ve^{2/3}A^2B-\ve B^3+\nonumber
\\
\big(\ve^{4/5}\t_0+
\ve^{5/6}T_1-\ve^{4/5}\sum_{n=1}^{\infty}\ve^{(2n-1)/5}\stackrel{n}{\th}_0
\big)B, \nonumber
\\
B'+3U_* A^2=-\ve^{1/3}(A^3-B^2U_*)-\ve^{2/3}AB^2+
\big(\ve^{1/6}T_1+\ve^{2/15}\t_0-
\nonumber\\
\ve^{2/15}\sum_{n=1}^\infty\ve^{(2n-1)/5}\stackrel{n}{\th}_0\big)
\big(U_*+\ve^{1/3}A\big),\ve\to0.
\label{3ineq}
\end{eqnarray}
\par
The asymptotics of $U$ as $\th\to-\infty$ rewriting in terms of
the variable $T_1$ gives the asymptotics of $A$ and $B$ as
$T_1\to-0$:
\begin{eqnarray}
A=\bigg({-2\over T_1^2}+\stackrel{2}{Y}\!\!^-_0
T_1^4+\dots\bigg)+
\ve^{2/15}\bigg({\t_0U_*^3\over5}T_1^2+\dots\bigg)+\nonumber\\
\ve^{1/6}\bigg({-1\over6}T_1^3+\dots\bigg)+\dots,\nonumber\\
B=\bigg({-2\over U_*^2T_1^3}+\stackrel{2}{Y}\!^-_0T_1^3+\dots\bigg)+
\ve^{2/15}\bigg({\t_0U_*\over5} T_1+\dots
\bigg)+\dots.\label{AB-as-left}
\end{eqnarray}
The dots in the brackets mean the terms of the less order with
respect to $T_1$ and out of the brackets -- with respect to
$\ve$.
\par
The system (\ref{3ineq}) and the asymptotics
(\ref{AB-as-left}) contain terms of the different orders with
respect to $\ve$. It is convenient to represent the formal
asymptotic solution with respect to $\ve$  in the form:
\begin{eqnarray}
A(T_1,\ve) = \sum_{n=0}^\infty\ve^{n/30}\stackrel{n}{A}(T_1), \nonumber\\
B(T_1,\ve) =
\sum_{n=0}^\infty\ve^{n/30}\stackrel{n}{B}(T_1),\quad \ve\to0.
\label{AB-anzatz}
\end{eqnarray}
\par
By substituting the representation for $A$ and $B$ in the
equation we obtain the recurrent sequence of the equations for
$\stackrel{n}{A}(T_1),\ \stackrel{n}{B}(T_1)$. It is
convenient for us to investigate the second order
differential equation for the functions
$\stackrel{n}{A}(T_1)$ instead of the system of differential
equations of the first order for
$\stackrel{n}{A}(T_1),\,\stackrel{n}{B}(T_1)$. The equation
for the leader term is

\bb
\stackrel{0}{A}\!\!^{\prime\prime} +
3\stackrel{0}{A}\!^2 = 0. \label{glp}
\ee
The equations for the next corrections of the asymptotics
(\ref{AB-anzatz}) are
$$
\stackrel{m}{A}\!''+6\stackrel{0}{A}\stackrel{m}{A}=\stackrel{m}{F},
$$
in particular,
\begin{eqnarray*}
\stackrel{m}{F}\equiv0,\,\,m=1,2,3,6,7;\quad
\quad\stackrel{4}{F}=-2U_*^2\t_0;\quad\stackrel{5}{F}=-2U_*^2T_1;\\
\stackrel{8}{F}=-3\stackrel{4}{A}\!\!^2(T_1);\quad\quad\quad
\stackrel{9}{F}=-6\stackrel{4}{A}(T_1)\stackrel{5}{A}(T_1);
\\ \stackrel{10}{F}=-\stackrel{1}{\th}_0
U_*-2U_*\bigg({\big(\stackrel{0}{A}\!\!'\big)^2\over2U_*}\,+
\,{3\stackrel{0}{A}\over2U_*^2}\bigg)-2U_*^2\bigg(\stackrel{0}{A}\!\!\!^3-
{\big(\stackrel{0}{A}\!\!'\big)^2\over4U_*^3}\bigg).
\end{eqnarray*}
\par
The asymptotics as $T_1\to-0$ of the coefficients of the formal
series for $A(T_1,\ve)$ can be obtained from the formulas
(\ref{AB-as-left})
\begin{eqnarray*}
\stackrel{0}{A}={-2\over T_1^2}+ \stackrel{2}{Y}\!\!^-_0
T_1^4+\dots;\quad
\stackrel{m}{A}\equiv0,\,\,m=1,2,3,6,7;\\
\stackrel{4}{A}={\t_0U_*^3\over5}T_1^2+\dots;\quad
\stackrel{5}{A}={-1\over6}T_1^3+\dots;\\
\stackrel{8}{A}(T_1)=O(T_1^6);\quad \stackrel{9}{A}(T_1)=O(T_1^7);
\quad \stackrel{10}{A}(T_1)=O(T_1^{-4}).
\end{eqnarray*}
The dots in the formulas mean the terms of the less order with
respect to $T_1$. In particular, these formulas give us the
domain where the segment of the length $n=10$ of the formal
series (\ref{AB-anzatz})  as $T_1\to-0$ is valid
$$
{\ve^{10/30}\stackrel{10}{A}(T_1)\over\stackrel{0}{A}(T_1)}\ll1,\,\,\,
{\hbox{or}}\,\,\, T_1\gg\ve^{1/6}.
$$
\par
Let us construct the formal asymptotic solution in the form
(\ref{AB-anzatz}). The seeking solution
$\stackrel{0}{A}(T_1)$ for the leader term equation
(\ref{glp}) has the asymptotics as $T_1\to -0$ in the form
$$
\stackrel{0}{A}=-{2\over T_1^2}+a_4(1) T_1^4+O(T_1^{10}),
$$
where $a_4(1)=\stackrel{2}{Y}\!\!^-_0$. Using this
asymptotics it is easy to see that the function
$\stackrel{0}{A}$ is expressed by the Weierstrass function
$$
\stackrel{0}{A}(T_1)=-2\wp(T_1;0,g_3(1)),\quad{\hbox{where}}\quad
g_3(1)=a_4(1)/56.
$$
\par
In a general case the solution $\stackrel{n}{A}(T_1)$ can be
represented in the neighborhood of the point $T_1=0$ in the
form
$$
\stackrel{n}{A}(T_1)=\stackrel{n}{A}_0\!\!\!^c(T_1) +
\stackrel{n}{x}_0\!\!\!^- A_1(T_1)+ \stackrel{n}{y}_0\!\!\!^-
A_2(T_1).
$$
The function  $\stackrel{n}{A}_0\!\!\!^c(T_1)$ is the partial
solution of the nonhomogeneous equation for the $n$-th correction.
This solution does not contain terms of the order of $T_1^{-3}$ and $T_1^4$.
The functions $A_1(T_1)$ and $A_2(T_1)$ are linear independent
solutions of the linear homogeneous equation in variations.
These solutions are uniquely defined by their asymptotics as
$T_1\to-0$
$$
 A_1(T_1)=T_1^{-3}+\dots,\quad
A_2(T_1)=T_1^4+\dots,
$$
the expansion of the function $A_1(T_1)$ does not contain
terms of the order of $T_1^4$. Both of these solutions can be represented
via  the Weierstrass function. In the case $g_2=0$
we obtain (\cite{Abramowitz})
$$
 A_1(T_1)={1\over4}\wp'(T_1,0,g_3),
$$ $$ A_2(T_1)=-14\pt_{g_3}\wp(T_1,0,g_3)\equiv
{7\over3g_3}T_1\wp'(T_1,0,g_3)+{14\over3g_3}\wp(T_1,0,g_3). $$
The Wronskian of these solutions is $ W=-{7\over4}$.
\par
The coefficients $\stackrel{n}{x}_0\!\!\!^-$ и
$\stackrel{n}{y}_0\!\!\!^-$ are defined by matching of the asymptotics
(\ref{AB-anzatz}) as $T_1\to-0$ and the asymptotics (\ref{exp-w}) as
$\th_0\to\infty$
$$
\stackrel{n}{x}_0\!\!\!^-=\stackrel{m}{X}_0\!\!\!^-,\quad n=6m+23;
$$ $$ \stackrel{n}{y}_0\!\!\!^-=\stackrel{m}{Y}_0\!\!\!^-,\quad
n=6m-12. $$
\par
The function $\stackrel{0}{A}$ has poles. One of them is
$T_1=0$. The second pole of this function is defined by an
elliptic integral
$$
\O_1=2\int_{A_{min}}^\infty{dy\over\sqrt{4y^3-g_3(1)}}, \quad
A_{min}=\bigg({g_3(1)\over4}\bigg)^{1/3}. $$
\par
The formal asymptotics (\ref{AB-anzatz}) loses the asymptotic
behaviour in the neighborhoods of poles. The asymptotics as
$T_1\to-0$ is known. These explicit formulas give us the
domain where the asymptotic series (\ref{AB-anzatz}) is valid
as $-T_1\gg\ve^{1/6}$.
\par
To construct the asymptotics of the coefficients of the formal
solution as $T_1\to-\O_1+0$ we solve the equations for the correction
terms and write their expansions as $T_1\to-\O_1+0$.
\par
It is convenient to represent the correction $\stackrel{n}{A}$ as
$T_1\to-\O_1+0$ in the form
\bb
\stackrel{n}{A}(T_1)=\stackrel{n}{A}_c(T_1)+
\stackrel{n}{x}_1\!\!\!^+ A_1(T_1+\O_1)+ \stackrel{n}{y}_1\!\!\!^+
A_2(T_1+\O_1). \label{n-correction}
\ee
The function $\stackrel{n}{A}_c(T_1)$ is the solution of the nonhomogeneous
equation for the $n$-th correction. The asymptotics as $T_1 \to -\O_1+0$ of
this function does not contain terms of the order of $(T_1+\O_1)^{-3}$ and
$(T_1+\O_1)^4$.
\par
The constants in the formula for the $n$-th correction are evaluated
by formulas
$$
\stackrel{n}{y}_1\!\!\!^+=\stackrel{n}{y}_0\!\!\!^-\,+\,\res_{r=0}\bigg[{1\over
r} \int_{-r}^{-\O_1+r} {dz\,\stackrel{n}{F}(z)A_1(z)\over W}
\bigg],
$$
$$
\stackrel{n}{x}_1\!\!\!^+= \stackrel{n}{x}_0\!\!\!^-\,+\,
C_1\stackrel{n}{y}_1\!\!^+\, +\,\res_{r=0}\bigg[{1\over r}
\int_{-r}^{-\O_1+r} {dz\,\stackrel{n}{F}(z)A_2(z)\over W}
\bigg],
$$
where the constant $C_1$ is defined by
$$
A_2(T_1+\O_1)=C_1A_1(T_1)+A_2(T_1).
$$
\par
To obtain the domain of the suitability of the formal
asymptotics as $T_1\to-\O_1+0$ it is necessary to know the
type of the singularity of the correction terms. The first
nonzero correction term is $\stackrel{4}{A}(T_1)$. The
expression (\ref{n-correction}) for $n=4$ contains the
partial solution $\stackrel{4}{A}_1\!\!\!^c(T_1)$. This
function is smooth and has the identical asymptotics as
$T_1\to-\O_1+0$ and as $T_1\to-0$. But the asymptotics of the
$\stackrel{4}{A}$ as $T_1\to-\O_1+0$ contains terms of the
order of $(T_1+\O_1)^4$ with the coefficients
$\stackrel{4}{x}_1\!\!\!^+\not=0$. Namely,
$$
\stackrel{4}{y}_1\!\!\!^+=0,
$$
$$
\stackrel{4}{x}_1\!\!\!^+={28\t_0\over3g_3}\z\big(\O_1/2;0,g_3\big),
$$
where $\z(T_1,g_2,g_3)$ is the Weierstrass $\z$-function.
\par
The values of the coefficients $\stackrel{n}{y}_1\!\!\!^+$ and
$\stackrel{n}{x}_1\!\!\!^+$ for $n>0$ also are evaluated. To
obtain the domain where the asymptotics
(\ref{AB-anzatz}) is valid it is necessary to note
$$
{\ve^{2/15}\stackrel{4}{A}(T_1)\over\stackrel{0}{A}(T_1)}=
O\left({\ve^{2/15}\over (T_1+\O_1)}\right)\ll 1.
$$
It is easy to see that the condition
$\ve^{-2/15}(\O_1+T_1)\gg1$ is the condition of the
suitability of the asymptotics (\ref{AB-anzatz}) as $T_1\to
-\O_1 + 0$.
\par
Thus, we have constructed the formal asymptotic solution of
the original equation in the form (\ref{AB-anzatz}) in the
domain where the following conditions $-\ve^{-1/6}T_1\gg1$
and $\ve^{-2/15}(\O_1+T_1)\gg1$ are valid.

\subsection{The bifurcation layer in the case of bounded  $k$}
\subsubsection{The neighborhoods of the poles $T_k$ of the
Weierstrass function}

In this section we construct the solution of the original
problem in the neighborhoods of the poles of the solutions of
(\ref{glp})  and show the change of the solution via the
passage from one pole to another one.
\par
In the neighborhood of the poles we will use a new scaling
variable
\bb
\th_k =
(T_k+\O_k-{1\over4}\sum_{n=1}^\infty\ve^{n/30}\stackrel{n}{x}_k\!\!\!^+)
\ve^{-1/6},\quad \ve\to0. \label{thk}
\ee
The original equation (\ref{sh}) looks like
\begin{eqnarray}
iU' + |U|^2U-\left(t_*+\ve^{4/5}\bigg(\tau_0\,+
\,\sum_{n=1}^\infty\ve^{(2n-1)/5}\stackrel{n}{\th}_0\bigg)
\,+ \, \right.  \nonumber
\\
\left. \ve^{5/6}\bigg(\sum_{j=1}^k\O_j \,+\,
{1\over4}\sum_{n=1}^{\infty}\ve^{n/30}
\sum_{j=1}^k\stackrel{n}{x}_j\!\!\!^+ \,+\, \ve^{1/6}
\th_k\bigg) \right)U=1,\quad \ve\to0, \label{Uk-eq}
\end{eqnarray}
where $\tau_0$ is the bifurcation point and  $\Omega_k$ are poles
of the Weierstrass function which describe the leader
term of the asymptotics in the intermediate domains.
\par
The solution is seeking in the form of formal series
\bb
U=U_* + \stackrel{0}{W}(\th_k) \,+\,
\ve^{4/5}\sum_{n=1}^\infty\ve^{(n-1)/30}\stackrel{n}{W}(\th_k),
\quad \ve\to0.
\label{weia}
\ee
\par
The coefficients of the formal expansion (\ref{weia}) are
satisfied to the recurrent system of equations. In
particular, the function $\stackrel{0}{W}(\th_k)$ is
satisfied to the equation (\ref{w0}). The equations for the
higher  corrections of the asymptotics are
\bb
i\stackrel{n}{W}{}^\prime +\left[2U_*\stackrel{0}{W}\!^* +
2U_*\stackrel{0}{W} - U_*^2 + 2 \stackrel{0}{W}
\stackrel{0}{W}\!^* \right]\stackrel{n}{W} + \left(U_* +
\stackrel{0}{W}\right)^2 \stackrel{n}{W}\!^* = \stackrel{n}{F}.
\label{enwaa}
\ee
\par
Let us become clear how the nonlinearity $\stackrel{n}{F}$
depends on $\stackrel{m}{W}$ and on
$\stackrel{j}{x}_k\!\!\!^+$. The maximum order of nonlinear
terms is cubic. The function $\stackrel{n}{F}$ contains terms
of the type
$$ \stackrel{j}{W}\stackrel{l}{W}\stackrel{m}{W}, \quad
n=j+l+m+46,\,\,j,l,m\ge1; $$
and
$$
(U_*+\stackrel{0}{W})\stackrel{j}{W}\stackrel{l}{W}, \quad
n=j+l+23, \,\,j,l\ge1; $$
linear terms with constant coefficients with dependness on
$\stackrel{j}{x}_l\!\!\!^+$
$$
\stackrel{m}{W}\sum_{l=1}^k \stackrel{j}{x}_l\!\!\!^+,\quad
n=m+j+25,l,\,m,j\ge1; $$
linear terms with the coefficients depending on $\th_k$
$$ \th_k\stackrel{m}{W},\quad n=m+55; $$
and with coefficients  only depending on
$\stackrel{m}{x}_j\!\!\!^+$
$$
\sum_{j=1}^k\stackrel{m}{x}_j\!\!\!^+(U_*+\stackrel{0}{W}),\quad
n=m+2; $$ and linear terms with depending on
$\stackrel{m}{\th}_0$:
$$
\stackrel{m}{\th}_0\stackrel{l}{W},\quad n=12m+l+18;
$$
and the terms depending on the leader term of the asymptotics
$$
\stackrel{m}{\th}_0(U_*+\stackrel{0}{W}),\quad n=12m-5.
$$
\par
For example, we write a few of the right hand-side terms for
the lower corrections:
$$
\stackrel{1}{F}=\tau_0 (U_* + \stackrel{0}{W});\quad
\stackrel{2}{F}=(U_* +
\stackrel{0}{W})\sum_{j=1}^k\O_j;\quad\,
\stackrel{n}{F}\equiv0,\,\,n=3,4,5,6; \quad
$$
$$
\stackrel{7}{F}=(U_* +
\stackrel{0}{W})(\th_k+{1\over4}\sum_{j=1}^k\stackrel{5}{x}_j\!\!\!^+).
$$
It is easy to see that the right-hand side of the equations
for a few of the lower corrections are proportional to leader
term of the asymptotics. The right-hand side for the 7-th
correction term depends on the leader term by more complicated
way. There is the variable coefficient $\th_k$ in the
equation.
\par
The solutions of the equations  for the corrections terms of
the asymptotics can be represented in the form
 \bb
\stackrel{n}{W}(\th_k)=\stackrel{n}{W}_c(\th_k) \,+\,
\stackrel{n}{X}_k\!\!\!^+\,w_1(\th_k) \,+\,
\stackrel{n}{Y}_k\!\!\!^+\,w_2(\th_k).\label{Wn} \ee Where
the function $\stackrel{n}{W}_c(\th_k)$ is the solution of
the nonhomogeneous equation for the $n$-th correction term.
This solution does not contain the terms of the order of
$\th_k^4$ and $\th_k^{-3}$ in its asymptotics as
$\th_k\to\infty$. The linear independent solutions of the
homogeneous equation $w_1(\theta_k)$ and $w_2(\theta_k)$ are
defined in the section 5.1.2. The constants
$\stackrel{n}{X}_k\!\!^+$ and $\stackrel{n}{Y}_k\!\!^+$ will
be determined below.
\par
The asymptotics of the function $\stackrel{0}{W}(\th_k)$ as
$\th_k\to\infty$ is given. It can be obtained by rewriting of the
formal solution $U_*+\ve^{2/6}A(T_k,\ve)+i\ve^{3/6}B(T_k,\ve)$ as
$T_k\to-\O_k+0$ in the terms of the scaling variable $\th_k$.
This asymptotics
 uniquely determines the leader term of the formal solution
 \bb
\stackrel{0}{W}(\th_k)={-2\over(\th_k-iU_*)^2}. \label{W0}
\ee
Thus,  the dependence on the number $k$ in the formula for
the leader term of the formal series is only contained in the
argument of the variable $\th_k$.
\par
To evaluate the constants $\stackrel{n}{X}_k\!\!\!^+$ and
$\stackrel{n}{Y}_k\!\!\!^+$ in the formula for the $n$-th
corrections we match its asymptotics as $\th_k\to\infty$ with
the asymptotics obtained by rewriting the formal series
$U_*+\ve^{2/6}A(T_k,\ve)+i\ve^{3/6}B(T_k,\ve)$ as
$T_k\to-\O_k+0$ in terms of the scaling variable $\th_k$. We
obtain
\bb
\stackrel{n}{X}_k\!\!\!^+\equiv0. \label{Xn}
\ee
The reason of this effect is that all terms connected
with the solutions of the type of $\stackrel{n}{x}_k\!\!\!^+A_1(T_k)$
have been matched by the shift of the independent
variable $\th_k$.
\par
To obtain the formula which connects the constants
$\stackrel{n}{Y}_k\!\!\!^+$ and $\stackrel{m}{y}_k\!\!\!^+$
it is necessary to note that solutions $A_2(T_k)$ is matching
with the solutions $w_2(\th_k)$. Thus the following formula
is valid
\bb
\stackrel{n}{Y}_k\!\!\!^+=\stackrel{m}{y}_k\!\!\!^+,\quad
n=m+7,\quad \stackrel{n}{Y}_k\!\!\!^+=0,\,n=0,1,\dots,6;
\label{Yn}
\ee
The value of the coefficient $\stackrel{7}{Y}_k\!\!\!^+$ is determined
by matching with the asymptotics of $\stackrel{0}{A}(T_k)$
\bb
\stackrel{7}{Y}_k\!\!\!^+=56 g_3(k). \label{Y6}
\ee
\par
Let us represent the expansions of two corrections of the (\ref{weia})
as $\th_k\to\infty$:
\bb
\stackrel{1}{W}(\th_k) = \stackrel{1}{W}_c(\th_k),\label{e1wasy}
\ee
where
$$
\stackrel{1}{W}_c(\th_k) = {{\tau_0 U_*^3}\over 5}\th_k^2 + i{{\tau_0
U_*}\over 5}
 \th_k - {\tau_0\over{30U_*}} + i{{3\tau_0}\over{5U_*^3}}\th_k^{-1} +
O(\th_k^{-2})
$$
These formulas explicitly contain the bifurcation point
$\tau_0$. The  terms of the asymptotics begin only  to depend
on the poles $\O_k,\ k=2,\dots$ of the Weierstrass
$\wp$-function with the second correction term
$\stackrel{2}{W}(\th_1)$. This correction as $\th_k\to\infty$
is represented in the form
\bb
\stackrel{2}{W}(\th_k) =\stackrel{2}{W}_c(\th_k) \label{e2wasy}
\ee
where
$$
\stackrel{2}{W}_c(\th_k)= {{P_k U_*^3}\over 5}\th_1^2 + i{{P_k
U_*}\over 5} \th_k - {P_k\over{30U_*}} +
i{{3P_k}\over{5U_*^3}}\th_k^{-1} + O(\th_k^{-2}),
$$
$$
P_k=\sum_{j=1}^k\O_j, \quad \th_\k \to \infty.
$$
\par
The formulas (\ref{thk}), (\ref{weia}), (\ref{Wn})--
(\ref{Y6}) determine the structure of the formal solution in
the neighborhood of the point $\th_k=+\infty$. It is easy to
see that  for $\forall k$ the formal solution (\ref{weia}) is
valid as $\th_k\ve^{1/6}\ll 1$.
\par
Now we know the form of the solution as $\th_k\to\infty$
and two linear independent solutions of the equation in variations
and we are able to write the representation of the correction
terms in the formal solution (\ref{weia}) as $\th_k\to-\infty$.
\bb
\stackrel{n}{W}(\th_k)=\stackrel{n}{W}_c(\th_k) \,+\,
\stackrel{n}{X}_k\!\!\!^-\,w_1(\th_k) \,+\,
\stackrel{n}{Y}_k\!\!\!^-\,w_2(\th_k). \label{Wn-} \ee
Where the functions $\stackrel{n}{W}_c(\th_k)$ have the similar
form as in the formula (\ref{Wn}).
The constants $\stackrel{n}{X}_k\!\!\!^-$ and $\stackrel{n}{Y}_k\!\!\!^-$
are
\bb
\stackrel{n}{X}_k\!\!\!^-=\stackrel{n}{\D}_k\!\!\!^X,\quad
\stackrel{n}{\D}_k\!\!\!^X=-\res_{r=\infty} \bigg[\,{1\over
r}\,\int_{-r}^{r}\big(\stackrel{n}{F}w_2^*+
\stackrel{n}{F}\!\!^*w_2\big){d\th_k\over{\cal W}}\bigg],
\label{Xn+}
\ee
\bb
\stackrel{n}{Y}_k\!\!\!^-=\stackrel{n}{Y}_k\!\!\!^+
+\stackrel{n}{\D}_k\!\!\!^Y,\quad \stackrel{n}{\D}_k\!\!\!^Y=-
\res_{r=\infty}\bigg[\,{1\over
r}\,\int_{-r}^{r}\big(\stackrel{n}{F}w_1^*+
\stackrel{n}{F}\!\!^*w_1\big){d\th_k\over{\cal W}}\bigg].
\label{Yn-}
\ee
It is important to note that
\bb
\stackrel{n}{Y}_k\!\!\!^-=0,
\quad{\hbox{as}}\quad n=1,\dots,6; \label{Y1-5} \ee \bb
\stackrel{7}{Y}_k\!\!\!^-=\stackrel{7}{Y}_k\!\!\!^+ +{\pi \over
2}. \label{Y0-}
\ee

\begin{demo} The formula (\ref{Y0-}) in the case of the
second order differential equation was obtained in the work
\cite{Dim-Hab}.
\end{demo}
\par
The domain where the formal series
(\ref{weia}) as $\th_k\to-\infty$ is valid
 is determined as well as
$\th_k\to\infty$. We obtain that for any value of $k$ this
domain is  determined by  $-\th_k\ve^{1/6} \ll 1$ as well as at
 large values $-\th_k$.

\subsubsection{The intermediate expansion, the internal
interval --- interval of the Weierstrass $\wp$-function}
\par
In this section we study the formal asymptotic solution in the
domain between $k$-th and $k+1$-th poles of the Weierstrass
function. The solution in this domain has the asymptotic
structure which is similar to  the asymptotic structure in
the initial intermediate interval. The difference is: the
asymptotics on the right side of the investigated interval is
matching with the asymptotics (\ref{weia}) where the
asymptotic sequence is $\{\ve^{n/30}\}$. But the formal
asymptotic solution has been constructed in the intermediate
domain on the initial interval and has been matched with the
asymptotic sequence  $\{\ve^{n/5}\}$.
\par
In this intermediate domain we use the slow variable
$T_k=\th_k\ve^{1/6}$ and seek the solution of the equation
(\ref{sh}) in the form

\bb
U=U_* + \ve^{2/6} A_k(T_k,\ve) + i \ve^{3/6}B_k(T_k,\ve).
\label{A-exp}
\ee
\par
The equation (\ref{sh}) becomes
\begin{eqnarray}
A_k'-2U_*^2 B_k=-\ve^{1/3}2U_*A_kB_k-\ve^{2/3}A_k^2B_k-\ve
B_k^3+\nonumber\\
\bigg(\ve^{4/5}\t_0-\ve^{4/5}\sum_{n=1}^{\infty}\ve^{(2n-1)/5}
\stackrel{n}{\th}_0
+\nonumber\\
\ve^{5/6}\big(T_k-\sum_{j=1}^k\O_j\big)\, -\,
\ve^{5/6}\sum_{n=1}^\infty\ve^{n/30}\sum_{j=1}^k
\stackrel{n}{x}_k\!\!\!^+ \bigg)B_k,
\nonumber\\
B_k'+3U_* A_k^2=-\ve^{1/3}(A_k^3-B_k^2U_*)-\ve^{2/3}A_kB_k^2+
\nonumber\\
\bigg(\ve^{2/15}\t_0+
\ve^{2/15}\sum_{n=1}^\infty\ve^{(2n-1)/5}\stackrel{n}{\th}_0\,+\,
\nonumber\\
 \ve^{1/6}\big(T_k-\sum_{j=1}^k\O_j\big)\, -\,
\ve^{1/6}\sum_{n=1}^\infty\ve^{n/30}\sum_{j=1}^k\stackrel{n}{x}_k\!\!\!^+
\bigg)\big(U_*+\ve^{1/3}A_k\big),\quad \ve\to0.
\label{3ineq-k}
\end{eqnarray}
\par
The asymptotics (\ref{weia}) as $\th_{k}\to-\infty$
of the function $U$ has been written in terms of the
new variable $T_k$ gives us the asymptotics as $T_k\to-0$ of the
functions $A$ and $B$
\begin{eqnarray}
A_k=\bigg({-2\over T_k^2}+\stackrel{2}{Y}\!^-_{k-1}
T_k^4+\dots\bigg)+
\ve^{2/15}\bigg({\t_0U_*^3\over5}T_k^2+\dots\bigg)+\nonumber
\\
\ve^{1/6}\bigg({P_kU_*^3\over5}T_k^2+{-1\over6}T_k^3+
\dots\bigg)+\dots,\nonumber\\
B_k=\bigg({-2\over U_*^2T_k^3}+\stackrel{2}{Y}\!^-_{k-1}
T_k^3+\dots\bigg)+ \ve^{2/15}\bigg({\t_0U_*\over5} T_k+\dots
\bigg)+\dots.\label{AB-as-left-k}
\end{eqnarray}
The dots in these formulas mean the terms of the small order
with respect to $T_k$  in the brackets and with respect to
$\ve$ out of them.
$$
P_k=\sum_{j=1}^k\O_j.
$$
\par
The system of the equations (\ref{3ineq-k}) and the asymptotics
(\ref{AB-as-left-k}) contain terms of the different orders with
respect on $\ve$ so we seek the formal asymptotic solution
in the form

\begin{eqnarray}
A_k(T_k,\ve) = \sum_{n=0}^\infty\ve^{n/30}\stackrel{n}{A}_k(T_k),
\nonumber\\
B_k(T_k,\ve) =
\sum_{n=0}^\infty\ve^{n/30}\stackrel{n}{B}_k(T_k),\quad
\ve\to0.
\label{AB-anzatz-k}
\end{eqnarray}
\par
Using the anzats for the functions $A_k$ and $B_k$ we obtain
the recurrent sequence of the equations for
$\stackrel{n}{A}_k(T_k),\ \stackrel{n}{B}_k(T_k)$. As it was
done above it is convenient for us to investigate the second
order differential equation for $\stackrel{n}{A}_k(T_k)$
instead of the system of the first order for
$\stackrel{n}{A}_k(T_k),\,\stackrel{n}{B}_k(T_k)$. The
equation for the leader term of asymptotics
(\ref{AB-anzatz-k}) is

\bb
\stackrel{0}{A}_k\!\!^{\prime\prime} + 3\stackrel{0}{A}_k\!^2 = 0.
\label{glp-k}
\ee
The equations for the next correction terms are

$$
\stackrel{m}{A}_k\!''+6\stackrel{0}{A}_k\stackrel{m}{A}_k=\stackrel{m}{F}.
$$
The expressions for the right hand terms $\stackrel{n}{F}$
depend on $n$. Here we represent the dependness of the
$\stackrel{n}{F}$ on corrections terms with a small
numbers and on the coefficients $\stackrel{j}{x}_k^+$.
The right hand functions $\stackrel{n}{F}$ contain
terms of the type
\newline
with quadratic dependence on corrections $\stackrel{m}{A}_k$
and $\stackrel{l}{A}_k$ as $n=m+l+10$;
\newline
with cubic dependence on corrections  $\stackrel{m}{A}_k$,
$\stackrel{j}{A}_k$  and  $\stackrel{l}{A}_k$ as $n=m+l+j+20$;
\newline
with linear dependence on corrections $\stackrel{m}{A}_k$ with
coefficients depending on $\stackrel{j}{x}_k^+$ as $n=m+j+15$;
\newline
with linear dependence on correction $\stackrel{m}{A}_k$ with
the linear coefficient with respect to
$(T_k+\sum_{j=1}^k\O_k)$ as $n=m+15$;
\newline
with linear dependence only on $\stackrel{j}{x}\!^+_k$ as
$n=j+5$.
\par
In particular,
\begin{eqnarray*}
\stackrel{m}{F}\equiv0,\,\,m=1,2,3,6,7;\quad
\quad\stackrel{4}{F}=-2U_*^2\t_0;\quad\stackrel{5}{F}=-2U_*^2T_k+P_k;\\
\stackrel{8}{F}=-3\stackrel{4}{A}_k\!\!^2(T_k);\quad\quad\quad
\stackrel{9}{F}=-6\stackrel{4}{A}_k(T_k)\stackrel{5}{A}(T_k);
\\ \stackrel{10}{F}=-\stackrel{1}{\th}_0
U_*-2U_*\bigg({\big(\stackrel{0}{A}_k\!\!'\big)^2\over2U_*}\,+
\,{3\stackrel{0}{A}_k\over2U_*^2}\bigg)
-2U_*^2\bigg(\stackrel{0}{A}_k\!\!\!^3-
{\big(\stackrel{0}{A}_k\!\!'\big)^2\over4U_*^3}\bigg).
\end{eqnarray*}
\par
The asymptotics as $T_k\to-0$ of the coefficients of the
formal series can be obtained from the formulas
(\ref{AB-as-left-k})
\begin{eqnarray*}
\stackrel{0}{A}_k={-2\over T_k^2}+ \stackrel{2}{Y}\!^-_k
T_k^4+\dots;\quad
\stackrel{m}{A}_k\equiv0,\,\,m=1,2,3,6,7;\\
\stackrel{4}{A}_k={\t_0U_*^3\over5}T_k^2+\dots;\quad
\stackrel{5}{A}_k={P_kU_*^3\over5}T_k^2+{-1\over6}T_k^3+\dots;\\
\stackrel{8}{A}_k(T_k)=O(T_k^6);\quad
\stackrel{9}{A}_k=O(T_k^6);\quad \stackrel{10}{A}=O(T_k^{-4}).
\end{eqnarray*}
In these formulas the dots mean the terms of the less order
with respect to $T_k$. In particular, from these formulas we
obtain the domain of suitability for the segment of the series
(\ref{AB-anzatz-k}) as $T_k\to-0$
$$
{\ve^{10/30}\stackrel{10}{A}_k(T_k)\over\stackrel{0}{A}_k(T_k)}\ll1,\,\,\,
{\hbox{or}}\,\,\, -\ve^{-1/6}T_k\gg1.
$$
\par
Let us construct the formal solution in the form
(\ref{AB-anzatz-k}). The seeking solution for the leader term
 of the asymptotics equation should have the asymptotic
structure  of the type as $T_k\to-0$
$$
\stackrel{0}{A}_k=-{2\over T_k^2}+\stackrel{2}{Y}\!^-_{k-1}
T_k^4+O(T_1^{10}).
$$
Using this formula it is easy to see that the leader term
$\stackrel{0}{A}_k$ of the asymptotics is expressed by the
Weierstrass $\wp$--function

$$
\stackrel{0}{A}_k(T_k)=-2\wp(T_k;0,g_3(k)),\quad{\hbox{where}}\quad
g_3(k)=\stackrel{2}{Y}\!^-_{k-1}/56.
$$
\par
In a general case the solution $\stackrel{n}{A}(T_k)$ can be
represented in the neighborhood of the point $T_k=0$ in the
form
$$
\stackrel{n}{A}_k(T_k)=\stackrel{n}{A}_k\!\!\!^c(T_k) +
\stackrel{n}{x}_k\!\!\!^- A_1(T_k)+ \stackrel{n}{y}_k\!\!\!^-
A_2(T_k).
$$
The function $\stackrel{n}{A}_k\!\!\!^c(T_k)$ is the partial
solution of the equation for the $n$-th correction term. This solution
does not contain terms of the order of  $T_k^{-3}$ and $T_k^4$ in the
asymptotics as $T_k\to 0$.
\par
The coefficients  $\stackrel{n}{x}_k\!\!\!^-$ and
$\stackrel{n}{y}_k\!\!\!^-$ are defined by matching of the
asymptotics (\ref{AB-anzatz-k}) as $T_k\to-0$ and the asymptotics
(\ref{weia}) as $\th_k\to\infty$
$$
\stackrel{n}{x}_k\!\!\!^-=\stackrel{m}{X}_k\!\!\!^-,\quad n=6m+23;
$$
$$
\stackrel{n}{y}_k\!\!\!^-=\stackrel{m}{Y}_k\!\!\!^-,\quad n=6m-12.
$$
\par
The period of the function  $\stackrel{0}{A}_k(T_k)$ is defined by
 the elliptic integral
$$
\O_k=2\int_{A_{min}}^\infty{dy\over\sqrt{4y^3-g_3(k)}}, \quad
A_{min}=\bigg({g_3(k)\over4}\bigg)^{1/3}.
$$
\par
The formal series solution  (\ref{AB-anzatz-k}) is not valid
in the neighborhood of the poles. The asymptotics as
$T_k\to0$ of the coefficients of the formal solution is
known.  These explicit formulas give us the domain where  the
formal solution (\ref{AB-anzatz-k}) is  valid.  This domain is
$-\ve^{-1/6}T_k\gg1$.
\par
To construct the asymptotics of the coefficients of the
formal series solution in the small neighborhood of the pole
$T_k=-\O_k$ we solve the equations for the corrections and
rewrite their expansions as $T_k\to-\O_k+0$
\par
It is convenient to represent the correction term
$\stackrel{n}{A}_k$ as $T_k\to-\O_k+0$ in the form

\bb
\stackrel{n}{A}_k(T_k)=\stackrel{n}{A}_c(T_k)+
\stackrel{n}{x}_{k+1}\!\!\!^+ A_1(T_k+\O_k)+
\stackrel{n}{y}_{k+1}\!\!\!^+ A_2(T_k+\O_k).
\label{n-correction-k}
\ee
Where the function $\stackrel{n}{A}_c(T_k)$ is the solution of
the nonhomogeneous equation for the $n$-th correction and its
expansion does not contain terms of the order of
$(T_k+\O_k)^{-3}$ and $(T_k+\O_k)^4$ as $T_k\to-\O_k$.
\par
The constants in the formula for the $n$-th corrections are
evaluated by
$$
\stackrel{n}{y}\!^+_{k+1}=\stackrel{n}{y}_k\!\!\!^-\,+\,
\stackrel{n}{\d}_k\!\!\!^y,\quad \stackrel{n}{\d}_k\!\!\!^y=
\res_{r=0}\bigg[{1\over r} \int_{-r}^{-\O_k+r}
{dz\,\stackrel{n}{F}(z)A_1(z)\over W} \bigg],
$$
$$
\stackrel{n}{x}\!^+_{k+1}= \stackrel{n}{x}_k\!\!\!^-\,+\,
C_k\stackrel{n}{y}\!^+_{k+1}\, +\, \stackrel{n}{\d}_k\!\!\!^x,
\quad \stackrel{n}{\d}_k\!\!\!^x=\res_{r=0}\bigg[{1\over r}
\int_{-r}^{-\O_k+r} {dz\,\stackrel{n}{F}(z)A_2(z)\over W}
\bigg],
$$
where the constant $C_k$ is defined by

$$
A_2(T_k+\O_k)=C_kA_1(T_k)+A_2(T_k).
$$
\par
To determine the domain of the suitability for the formal
series (\ref{AB-anzatz-k}) as $T_k\to-\O_k+0$ it is important
to know the structure and orders of the singularities of the
correction terms. The first nonzero correction term is
$\stackrel{4}{A}(T_k)$. The representation
(\ref{n-correction-k}) for $n=4$ contains the function
$\stackrel{4}{A}_{k+1}\!\!\!^c(T_k)$  as the partial solution
of the correction equation. This solution has the identical
structure of the asymptotics as $T_1\to\O_1+0$ and as
$T_k\to-0$. But the coefficient
$\stackrel{4}{x}_1\!\!\!^+\not=0$ as $T_1\to\O_1+0$. Namely,

$$
\stackrel{4}{y}_{k+1}\!\!\!^+=0,
$$
$$
\stackrel{4}{x}_{k+1}\!\!\!^+={28\t_0\over3g_3}\z\big(\O_k/2;0,g_3\big),
$$
where $\z(T_1,g_2,g_3)$ is the Weierstrass $\z$-function.
\par
The values of the coefficients $\stackrel{n}{y}_1\!\!\!^+$ and
$\stackrel{n}{x}_1\!\!\!^+$ for $n>0$ can be also evaluated.
To obtain the domain of the suitability of the formal series
(\ref{AB-anzatz-k}) we note
$$
{\ve^{2/15}\stackrel{4}{A}(T_k)\over\stackrel{0}{A}(T_k)}=
{\ve^{2/15}\over (T_k+\O_k)}\const\ll1.
$$
It is easy to see that the condition
$\ve^{-2/15}(T_k+\O_k)\gg1$ describes the domain where this
solution (\ref{AB-anzatz-k}) is valid as $T_k\to-\O_k$.
\par
Thus in this section we have constructed the formal series
solution in the form (\ref{AB-anzatz-k}) and have shown that
this representation of the solution is valid as
$-\ve^{-1/6}T_k\gg1$ and $\ve^{-2/15}(T_k+\O_k)\gg1$.
\par

\subsubsection{The discrete dynamical system}

\par
The coefficients which determine the intermediate and the
second internal solution are changed by passage from the $k$
to $k+1$. Let us investigate the discrete dynamical system for
the coefficients $\stackrel{n}{X}_k\!\!\!^\pm$,
$\stackrel{n}{Y}_k\!\!\!^\pm$ and $\th_k$,
$\stackrel{n}{x}_k\!\!\!^-$, $\stackrel{n}{y}_k\!\!\!^\pm$.
\par
The connection formulas  for these coefficients are
$$
\stackrel{n}{y}\!^+_{k+1}=
\stackrel{m}{Y}\!^-_k\,+\,\stackrel{n}{\delta}\!^y_k,
$$
where $n=6m-12$,
$$
\stackrel{m}{Y}\!^-_{k+1}=\stackrel{n}{y}_k\!\!\!^+\,+\,
\stackrel{m}{\Delta}\!^Y_{k+1}.
$$
Where $\stackrel{m}{Y}_k\equiv0$ as $n=0,\dots,6$. In the
case of $m=7$ we determine the parameter $g_3(k)$ of the
Weierstrass $\wp$--function for the leader term of the
intermediate expansion
$$
g_3(k)={1\over56}\big(a_4+ {\pi\over2} (k-1)\big).
$$
The connection formula is defined by explicit way without
investigating of the next corrections because
at first
$$
\stackrel{0}{\delta}_k\!\!\!^y=0.
$$
and secondly
$$
\stackrel{7}{\Delta}\!^Y_{k+1}={\pi\over2}.
$$
\begin{demo} This connection formula for $g_3(k)$  in the case
of second order ordinary differential equation  was found in
the work \cite{Dim-Hab} without investigation of high
correction terms in internal domains and discrete dynamical
system.
\end{demo}
\par
It is easy to see that this system for
$\stackrel{m}{Y}\!^-_{k+1}$ and
$\stackrel{n}{y}_{k+1}\!\!\!^+$ is not closed. It is
connected with the fact that the terms
$\stackrel{m}{\Delta}\!^Y_{k+1}$ and
$\stackrel{n}{\delta}\!^y_{k+1}$ contain the terms
$\stackrel{l}{X}_k\!\!\!^-$ and $\stackrel{j}{x}_k\!\!\!^\pm$.
So to obtain the closed system for the coefficients two above
relations have to supplement by equations for
$\stackrel{l}{X}_k\!\!\!^-$ and $\stackrel{j}{x}_k\!\!\!^\pm$.
As result we obtain
$$
\stackrel{n}{x}_k\!\!\!^+=C_k\stackrel{n}{y}_k\!\!\!^+\,+\,
\stackrel{n}{\delta}{}^x_{k-1}
\,+\,\stackrel{n}{\Delta}{}^X_{k-1}.
$$
\par
This discrete dynamical system allows us to evaluate the
phase variable in the leader term of the second internal
expansion
$$
\th_{k+1}=\th_k+\ve^{-1/6}\bigg(\O_{k+1}+
{1\over4}\sum_{n=1}^\infty
\ve^{n/30}\stackrel{n}{x}_{k+1}^+\bigg),\quad  \ve\to0.
$$
\par
Thus, to obtain the phase variable to within $o(1)$ in the
neighborhood of the  $(k+1)$-th pole it is necessary to
obtain the term $\stackrel{n}{x}_{k+1}\!\!\!^+$ as $n\le5$ and
$\stackrel{m}{Y}_k\!\!\!^-$ as $m\le3$.
\par
Let us  consider these formulas as $n\le5$ and $m\le3$
$$
\stackrel{m}{\D}_{k}\!\!\!^Y\equiv0,\quad m\not=7; $$ $$
\stackrel{n}{\d}_k\!\!\!^y\equiv0,\quad 0<n<5; $$ $$
\stackrel{5}{\d}_k\!\!\!^y=2\z(\O_k/2,0,g_3(k)). $$
We also have that
 $\stackrel{l}{\D}_k\!\!\!^X\equiv0,$ as $l<1$.
\par
As result the formulas for $\stackrel{n}{x}_k\!\!\!^+$
can be written in the form
$$
 \stackrel{n}{y}_k\!\!\!^+\equiv0\quad
\stackrel{n}{x}_k\!\!\!^+\equiv0,\quad n=1,2,3,$$ $$
\stackrel{4}{y}_k\!\!\!^+\equiv0,\quad
\stackrel{4}{x}_k\!\!\!^+=\stackrel{4}{\d}_k\!\!\!^x,$$ $$
\stackrel{5}{y}_k\!\!\!^+=\stackrel{5}{y}\!^+_{k-1} \,+\,
\stackrel{5}{\d}\!^y_{k-1},\quad
\stackrel{5}{x}\!^+_k=C_k\stackrel{5}{y}\!^+_k \,+\,
\stackrel{5}{\d}\!^x_k.
$$

\subsubsection{The domain of the suitability of the intermediate
expansion}

\par
The domain where the intermediate asymptotic
expansion is valid for big values of $k$ is defined by inequality
 $$
{\ve^{1/6}\stackrel{5}{A}_k\over\stackrel{0}{A}_k}\ll1.$$
\par
The function $\stackrel{5}{A}_k$ contains the term
$\stackrel{5}{y}_k\!\!\!^+ A_1(T_k)$ which is the biggest
with respect to $k$
$$
\stackrel{5}{y}_k\!\!\!^+= \sum_{j=1}^k\stackrel{5}{\d}_j\!\!\!^y
=O(k^{7/6})$$
\par
This formula leads us to the inequality
 $$ \ve^{1/6}k^{7/6}\ll1,\quad
k\ll\ve^{-1/7}.
$$
It is easy to determine the order of
$P_k$ as $k\to \infty$
$$
P_k=\sum_{j=1}^k\O_j=O(k^{5/6}),\quad \hbox{при}\quad k\to\infty.
$$

\subsection{The intermediate expansion for large $k$}

\par
In this section to take into account the secular terms as
large values of $k$ we modify the intermediate expansion.

\subsubsection{The domain for moderate large values of $k$}

\par
As was shown above for the large values of $k$ the correction
terms of the asymptotics contain the secular terms with
respect to $k$. These singularities have the different
structure. At first, the singularities are connected with the
terms $\stackrel{n}{y}_k\!\!\!^+$ as the solutions of the
homogeneous equations. Secondly, the singular terms $P_k$ and
$\sum_{j=1}^k\stackrel{n}{x}_j\!\!\!^+$ are contained in
right-hand sides of the equations (\ref{3ineq-k}). These
singularities lead to lose of the suitability of the
equations for the correction terms equations.
\par
To suppress these singularities we allow the dependness
of the parameter $g_3$ on the small parameter
$$
g_3(k,\ve)=\stackrel{0}{g}_3(k)+\sum_{n=1}^\infty
\ve^{n/30}\stackrel{n}{g}_3(k).
$$
Where
$$
\stackrel{n}{g}_3(k+1)=\stackrel{n}{g}_3(k)+{1\over56}\big(
\stackrel{m}{\D}_k\!\!\!^Y+\stackrel{n}{\d}_k\!\!\!^y\big),\quad m=n+5.
$$
\par
The equations for $\stackrel{n}{x}_k\!\!\!^+$ become
$$
\stackrel{n}{x}_k\!\!\!^+=C_k\stackrel{n}{\d}_k\!\!\!^y \,+\,
\stackrel{l}{\D}\!^X_k,\quad n=l+29.
$$
We can suppress the first kind of the singularities in the
intermediate expansion by this way. To suppress the
singularities of the second type connected with large values
of $P_k$ and sum of the coefficients
$\stackrel{n}{x}_j\!\!\!^+$ we include the following term in
equation for the leader term of the intermediate expansion
$$
\l_k(\ve)=\ve^{1/6}\bigg(P_k+{1\over4}\sum_{n=1}^\infty\ve^{(n-1)/30}
\sum_{j=1}^k\stackrel{n}{x}_j\!\!\!^+\bigg).
$$
Then the leader term equation has the form
$$
\stackrel{0}{A}_k''+3\stackrel{0}{A}_k\!\!\!^2=\l_k,
$$
\par
and the solution is
$$
\stackrel{0}{A}_k(T_k)=-2\wp(T_k,\l_k/2,g_3(k,\ve)).
$$
\par
Including of the term $\l_k(\ve)$ in equation for the leader term
of the asymptotics leads us to change of the asymptotic sequence.
Namely, $$ A(T_k,\ve)=\sum_{n=0}^\infty
\ve^{n/6}\stackrel{5n}{A}(T_k). $$ The terms of the order of
$\ve^{m/30}$ with $m\not= 5l, l\in Z$ are contained in the leader
term of the asymptotics and corrections by introducing of the
$\l_k(\ve)$. To obtain the intermediate expansion with respect to
$\ve^{m/30}$ one should expand this representation with respect
to $\ve$.
\par
The solutions of the corrections equations are defined by
their asymptotics as $T_k\to-0$:
$$
\stackrel{5n}{A}_k(T_k)=\stackrel{5n}{A}_k\!\!\!^c(T_k)+
\stackrel{5n}{x}_k\!\!\!^- A_1(T_k)
$$
The main term of the asymptotics with respect to $T_k$ is the
term of the order of $T_k^4$. The function
$\stackrel{n}{A}_k$ is the partial solution of the
nonhomogeneous equation for the correction term and does not
contain the terms of the order of $T_k^4$ and $T_k^{-3}$ in
the asymptotics as $T_k\to-0$. The function $A_1$ is the
solution of variation equation and is defined by the leader
term
$$
A_1={d\over dT_k}\wp(T_k,\l_k/2,g_3(k,\ve)).
$$
\par
In the neighborhood of the point $T_k=-\O_k+0$ the corrections
can be represented by formulas
$$
\stackrel{n}{A}_k(T_k)=\stackrel{n}{A}_k\!\!\!^c(T_k)+
\stackrel{n}{x}\!^+_{k+1} A_1(T_k) \,+\,
\stackrel{n}{\d}\!^y_k A_2(T_k).
$$
The function $A_2$ is the second solution of the homogeneous
variation equation
$$
A_2={d\over dg_3}\wp(T_k,\l_k/2,g_3),\quad \hbox{при}\quad
g_3=g_3(k,\ve).
$$
\par
The constructed solution is valid between the poles of the
Weierstrass $\wp$-function $\wp(T_k,\l_k/2,g_3(k,\ve))$ in
the similar domains as in the case of intermediate expansion
for finite values of $k$.

\subsubsection{The second internal expansion for moderate
large $k$}
\par
To construct the second internal expansion in this case
we  should make the changes in the formula for evaluating of
 $\stackrel{n}{Y}_k\!\!\!^+$:
$$
\stackrel{n}{Y}_k\!\!\!^+=56\stackrel{m}{g}_3(k), \quad n=6m-12.
$$

\subsubsection{The large values of $k$}

\par
The large values of $k$ lead to grow of the term $\l_k(\ve)$
in the  right hand side of the equation for leader term of the
asymptotics. The structure of the formal solution in this
case becomes complicated. Two parameters appear and the
formal serires solution of the original equation use both of
them. These two parameters are  $\ve$ as small one and
$\l_k(\ve)$ as  large one. To simplify the representation of
the formal series solution in this section we only use the
segments of the formal asymptotic series. We describe the
domain of the suitability  of this representation of the
formal series solution. It is convenient to construct the
leader term of the intermediate asymptotics in the form
$$
A_k=\sqrt{\l_k}\big(\stackrel{0}{\cal A}_k \,+\,
\ve^{1/6}\stackrel{5}{\cal A}_k\big) \,+\,
\ve^{1/3}\l_k\big(\stackrel{10}{\cal A}_k \,+\,
\ve^{1/2}\stackrel{15}{\cal A}_k\big) \,+\,
\ve^{2/3}\l_k^{3/2}\stackrel{20}{\cal A}_k,
$$
$$
B_k=\l_k^{3/4}\big(\stackrel{0}{\cal B}_k \,+\,
\ve^{1/6}\stackrel{5}{\cal B}_k\big) \,+\,
\ve^{1/3}\l_k^{5/4}\big(\stackrel{10}{\cal B}_k \,+\,
\ve^{1/6}\stackrel{15}{\cal B}_k\big) \,+\,
$$
$$
\ve^{2/3}\l_k^{7/4}\big(\stackrel{20}{\cal B}_k+
\ve^{1/6}\stackrel{25}{\cal B}_k\big).
$$
The coefficients of the asymptotics contain the independent
variable with the large term $\s_k=\l_k^{1/4}T_k$, namely,
$\stackrel{n}{\cal A}_k=\stackrel{n}{\cal A}_k(\s_k,k,\ve)$
and
 $\stackrel{n}{\cal B}_k=\stackrel{n}{\cal B}_k(\s_k,k,\ve)$.
\par
The leader term of the asymptotics satisfies to the equation
\begin{eqnarray*}
\stackrel{0}{\cal A}_k'-2U_*\stackrel{0}{\cal B}_k=0,\\
\stackrel{0}{\cal B}_k'+3U_*\stackrel{0}{\cal A}_k\!\!\!^2= U_*.
\end{eqnarray*}
It is easy to write the equations for the correction terms
of the asymptotics
\begin{eqnarray*}
\stackrel{5}{\cal A}_k'-2U_*\stackrel{5}{\cal B}_k=0,\\
\stackrel{5}{\cal B}_k'+6U_*\stackrel{0}{\cal A}_k\stackrel{5}{\cal A}_k=
U_*T_k/\l_k;
\end{eqnarray*}
\begin{eqnarray*}
\stackrel{10}{\cal A}_k'-2U_*\stackrel{10}{\cal B}_k=
-2U_*\stackrel{0}{\cal A}_k\stackrel{0}{\cal B}_k,\\
\stackrel{10}{\cal B}_k'+6U_*\stackrel{0}{\cal
A}_k\stackrel{10}{\cal A}_k= -3U_*\stackrel{5}{\cal A}_k\!\!\!^2
\,+\, \stackrel{0}{\cal A}_k - \stackrel{0}{\cal A}_k\!\!\!^3
-U_*\stackrel{0}{\cal B}_k\!\!\!^2;
\end{eqnarray*}
\begin{eqnarray*}
\stackrel{15}{\cal A}_k'-2U_*\stackrel{15}{\cal B}_k=
-2U_*\big(\stackrel{5}{\cal A}_k\stackrel{0}{\cal B}_k+
\stackrel{0}{\cal A}_k\stackrel{5}{\cal B}_k\big),\\
\stackrel{15}{\cal B}_k'+6U_*\stackrel{0}{\cal A}_k\stackrel{15}{\cal A}_k=
\big(1-3\stackrel{0}{\cal A}_k\!\!\!^2\big)\stackrel{5}{\cal A}_k
\,+\,
{T_k\over\l_k}\stackrel{0}{\cal A}_k
\,- \,
2U_*\stackrel{0}{\cal B}_k\stackrel{5}{\cal B}_k
\,-\,
6U_*\stackrel{5}{\cal A}_k\stackrel{10}{\cal A}_k;
\end{eqnarray*}
\begin{eqnarray*}
\stackrel{20}{\cal A}_k\!'-2U_*\stackrel{20}{\cal B}_k=
-\stackrel{0}{\cal B}_k-\stackrel{0}{\cal A}_k\!\!\!^2
\stackrel{0}{\cal B}_k -
\\
2U_*\big(\stackrel{10}{\cal A}_k\stackrel{0}{\cal B}_k \,+\,
\stackrel{0}{\cal A}_k\stackrel{10}{\cal B}_k \,+\,
{1\over\sqrt{\l_k}}\stackrel{5}{\cal A}_k\stackrel{5}{\cal
B}_k\big),
\\
\stackrel{20}{\cal B}_k\!'+6U_*\stackrel{0}{\cal
A}_k\stackrel{20}{\cal A}_k= \stackrel{10}{\cal
A}_k-3\stackrel{0}{\cal A}_k\!\!\!^2\stackrel{10}{\cal A}_k-
\stackrel{0}{\cal A}_k\stackrel{0}{\cal B}_k\!\!\!^2-
\\
2U_*\stackrel{0}{\cal B}_k\stackrel{10}{\cal B}_k \,+\,
{1\over\l_k^{1/2}} \big(-3 \stackrel{0}{\cal A} \stackrel{5}{\cal
A}_k\!\!\!^2 - U_*\stackrel{5}{\cal B}_k\!\!\!^2\big)
\\
\,+\,{1\over\l_k^{3/2}} T_k\stackrel{5}{\cal A}_k
\,-\,
3U_*\stackrel{10}{\cal A}_k\!\!\!^2- 6U_*\stackrel{5}{\cal
A}_k\stackrel{15}{\cal A}_k.
\end{eqnarray*}
\par
As was shown above it is convenient to investigate the second order
differential equation for the functions $\stackrel{n}{\cal A}_k$.
In particular, we have

$$
\stackrel{0}{\cal A}_k\!\!\!''+3\stackrel{0}{\cal A}_k\!\!\!^2=1.
$$
The solution of this equation matched with the asymptotics of the internal
formal solution in the neighborhood of the $k-1$-th pole is
$$
\stackrel{0}{\cal A}=-2\wp(\s_k,1/2,\g_3(k,\ve)),
\quad\hbox{where}\quad \g_3(k,\ve)={g_3(k,\ve)\over\l_k^{3/2}}.
$$
\par
The parameter $g_3(k,\ve)$ and corrections of the asymptotics are
obtained as was shown above in the domain of the intermediate
expansion. It is easy to see that unique type of singularities
is connected with the parameter $\l_k$. The intermediate
expansion is valid as
$$
{\ve^{1/3}\l_k\over\sqrt{\l_k}}\ll1,\quad \hbox{or}\quad \l_k\ll\ve^{-2/3}.
$$
\par
The constructed sequence of the intermediate and the second internal
expansions is valid in the domain as
$$
(t_*-t)\ll\ve^{1/6}.
$$

\section{Fast oscillating asymptotic solution}

In this section we construct a special fast oscillating
asymptotic solution of the equation (\ref{sh}) by
Krylov-Bogolyubov method \cite{K-B}, \cite{Kuz}. This
constructed solution is valid as $(t_*-t)\ve^{-2/3} \gg 1$.
This special solution degenerates at the point $t=t_*$. In the
neighborhood of this point this special solution loses the
fast oscillating structure and matchs with the asymptotics
from the previous section.

\subsection{The family of the fast oscillating solutions}

Here we construct the family of the fast oscillating solutions
which depends on two parameters.
\par
The fast oscillating solution of the equation
(\ref{sh}) is seeking in the form
\bb
U(t,\ve) = \stackrel{0}{U}(t_1,t,\ve) + \ve
\stackrel{1}{U}(t_1,t,\ve) + \ve^2 \stackrel{2}{U}(t_1,t,\ve) +
\dots, \label{as4}
\ee
where $t_1=S(t)/\ve+\phi(t)$ is the fast variable, the functions
$S(t)$ and $\phi(t)$ are unknown functions.
\par
The problem on constructing the fast oscillating solution for
the second order differential equation was investigated by
many authors. In this section we construct the solution in
the  manner of F.J.Bourland and R.Haberman \cite{Bour-Hab2}.
They have investigated the second order equation so the
direct reference for the equation (\ref{sh}) is not valid. In
this section we represent some evaluations which lead us  to
the elegant formulas from the work \cite{Bour-Hab2}.
\par
Let us substitute (\ref{as4}) in equation (\ref{sh}). As a result
we obtain the equations for the leader term of the asymptotics
(\ref{as4}) and the first correction term. \bb
iS^{\prime}\pt_{t_1} \stackrel{0}{U} +
|\stackrel{0}{U}|^2\stackrel{0}{U} - t \stackrel{0}{U} = 1.
\label{u0} \ee \bb iS^{\prime}\pt_{t_1} \stackrel{1}{U} +
\left(2|\stackrel{0}{U}|^2-t\right) \stackrel{1}{U} +
\stackrel{0}{U}\!^2\stackrel{1}{U^*} = -i\pt_t \stackrel{0}{U}
-i\phi'\stackrel{0}{U}. \label{u1} \ee
\begin{eqnarray}
iS^{\prime}\pt_{t_1} \stackrel{2}{U} +
\left(2|\stackrel{0}{U}|^2-t\right) \stackrel{2}{U} +
\stackrel{0}{U}\!^2\stackrel{2}{U^*} = \nonumber
\\
 -i\pt_t
\stackrel{1}{U}-i\phi'\pt_{t_1}\stackrel{1}{U}-
2|\stackrel{1}{U}|^2\stackrel{0}{U}-
(\stackrel{1}{U})^2\stackrel{0}{U}\!^*. \label{u2}
\end{eqnarray}
\par
The time variable $t$ is the independent variable in these
equations. The equation (\ref{u0})  has the first integral
with respect to "fast" variable $t_1$

\bb
{1\over 2}|\stackrel{0}{U}|^4 -t |\stackrel{0}{U}|^2 -
(\stackrel{0}{U} + \stackrel{0}{U^*}) = E(t), \label{int}
\ee

where $E(t)$ is a "constant" of integration.

\par
This expression for the first integral can be considered as
the equation for the function $\stackrel{0}{U}$. Let us
express $\stackrel{0}{U}$ through the complex conjugate
function and as a result we obtain the equation only for
$\stackrel{0}{U}$.

$$
iS^{\prime}\pt_{t_1} \stackrel{0}{U} = \pm
\sqrt{2\stackrel{0}{U}\!^3 + (2E(t)+t^2) \stackrel{0}{U}\!^2 +
2t\stackrel{0}{U}+1}.
$$
This equation can be easy integrated

\bb
iS'\int_{u_0}^{\stackrel{0}{U}}
{dy\over\pm\sqrt{2y^3+(2E+t^2)y^2+2ty+1}}= t_1+S_0.
\label{u0s}
\ee
Where $u_0$ and $S_0$ are constants and we integrate on the
curve $\G(t)$ on the complex plane $y$. This curve is defined
by

$$
{1\over 2}|y|^4 -t |y|^2 - (y + y^*) = E(t).
$$
The sign $+$ or  $-$ fixed the sheet of a Riemann surface
where the initial point $u_0$ is situated. We choose the "$+$"
sign.

\par
The motion defined by (\ref{u0s}) is periodic with respect to
fast variable $t_1$. By integrating of the formula (\ref{u0s})
over the period we obtain:
\bb
iS'\int_{\G(t)}{dy\over\sqrt{2y^3+(2E+t^2)y^2+2ty+1}}=T.
\label{period}
\ee
The constant $T$ is the period of the function $\stackrel{0}{U}$
with respect to fast variable $t_1$. This formula
is the differential equation for the function $S(t)$.

\par
It is convenient to represent the leader term
$\stackrel{0}{U}(t_1,t)$ as the sum of the real and the
imaginary parts
$$
\stackrel{0}{U}=\stackrel{0}{U}_R+i\stackrel{0}{U}_I.
$$
By shift of phase variable $t_1$ we can obtain that
$\stackrel{0}{U}_R(t_1,t)$ is even function with respect
to $t_1$ and $\stackrel{0}{U}_I(t_1,t)$ is odd one.
This property of the functions $\stackrel{0}{U}_I(t_1,t)$
and $\stackrel{0}{U}_R(t_1,t)$ is convenient
to evaluate the integrals with respect to fast
variable.

\par
To determine the leader term of the asymptotics (\ref{as4})
it is necessary to obtain an equation for $E(t)$. This
equation appears from a condition of the boundedness of the
first correction term $\stackrel{1}{U}$ of the asymptotics.
\par
The general solution of the equation for the first correction
can be represented by two linear independent solutions of the
homogeneous equation

$$ u_1(t_1,t)=\pt_{t_1}\stackrel{0}{U},\quad \hbox{and}\quad
u_2(t_1,t)={\pt_ES\over S'}t_1 u_1 +\pt_E \stackrel{0}{U}. $$ In
these formulas the $u_1$ is periodic bounded function, but $u_2$
grows linearly with respect to $t_1$. The Wronskian of these
solutions is

$$
W=u_1u_2^*-u_1^*u_2={1\over iS'}.
$$
It can be evaluated by explicit formulas
$$
W={1\over iS'} \big(-(|\stackrel{0}{U}|^2-t)\stackrel{0}{U}+1\big)
\pt_E\stackrel{0}{U}\!^* + {1\over iS'}
\big(-(|\stackrel{0}{U}|^2-t)\stackrel{0}{U}\!^*+1\big)
\pt_E\stackrel{0}{U}=
$$
$$
= {1\over iS'}\pt_E\big(-{1\over2}
|\stackrel{0}{U}|^4+t|\stackrel{0}{U}|^2+\stackrel{0}{U}
+\stackrel{0}{U}\!^*\big)={-1\over iS'}.
$$
\par
To obtain necessary condition of the boundedness   of the
first correction terms of the asymptotics (\ref{as4}) we
multiply the equation (\ref{u1}) to
$\pt_{t_1}\stackrel{0}{U}\!^*$ and integrate it with respect
to the variable $t_1$ over the  period. Then we integrate the
first term by parts

$$
iS'\int_{0}^Td t_1 \pt_{t_1}\stackrel{1}{U}\pt_{t_1}
\stackrel{0}{U}\!^*=-\int_0^T dt_1\bigg(2(|\stackrel{0}{U}|^2-t)
\pt_{t_1}\stackrel{0}{U}\!^*+\big(\stackrel{0}{U}\!^*\big)^2
\pt_{t_1}\stackrel{0}{U}\bigg)\stackrel{1}{U}.
$$
As a result we obtain:

$$
\int_0^Tdt_1\bigg[\stackrel{0}{U}\!^2\stackrel{1}{U}\!^*
\pt_{t_1}\stackrel{0}{U}\!^*-\big(\stackrel{0}{U}\!^*\big)^2
\pt_{t_1}\stackrel{0}{U}\stackrel{1}{U}\bigg]=
$$
$$
-i\int_0^Tdt_1\pt_t\stackrel{0}{U}\pt_{t_1}\stackrel{0}{U}\!^*
-i\phi'\int_0^Tdt_1
\pt_{t_1}\stackrel{0}{U}\pt_{t_1}\stackrel{0}{U}\!^*.
$$
\par
Let us  consider the equation for the complex conjugate
function

$$
-iS^{\prime}\pt_{t_1} \stackrel{1}{U}\!^* +
\left(2|\stackrel{0}{U}|^2-t\right) \stackrel{1}{U}\!^* +
\big(\stackrel{0}{U}\!^*\big)\!^2\stackrel{1}{U}= i\pt_t
\stackrel{0}{U}\!^*+i\phi\pt_{t_1}\stackrel{0}{U}\!^*.
$$
Multiply this equation to $\pt_{t_1}\stackrel{0}{U}$ and
integrate it with respect to $t_1$ over the  period. Then we
obtain
$$
\int_0^Tdt_1\bigg[\big(\stackrel{0}{U}\!^*\big)^2\stackrel{1}{U}
\pt_{t_1}\stackrel{0}{U}-\stackrel{0}{U}\!^2
\pt_{t_1}\stackrel{0}{U}\!^*\stackrel{1}{U}\!^*\bigg]=
$$
$$
i\int_0^Tdt_1\pt_t\stackrel{0}{U}\!^*\pt_{t_1}\stackrel{0}{U} +
i\phi'\int_0^Tdt_1
\pt_{t_1}\stackrel{0}{U}\pt_{t_1}\stackrel{0}{U}\!^*.
$$
Combining obtained expressions we get

$$
-i\int_0^Tdt_1 \pt_t\stackrel{0}{U}\pt_{t_1}\stackrel{0}{U}\!^*+
i\int_0^Tdt_1 \pt_t\stackrel{0}{U}\!^*\pt_{t_1}\stackrel{0}{U}=0.
$$
After integrating by parts we get

$$
i\int_0^Tdt_1 \stackrel{0}{U}\!^*\pt_{t_1}\pt_t\stackrel{0}{U}+
i\int_0^Tdt_1\pt_t\stackrel{0}{U}\!^*\pt_{t_1}\stackrel{0}{U}=0,
$$
or
$$
i\pt_t\int_0^Tdt_1\stackrel{0}{U}\!^*\pt_{t_1}\stackrel{0}{U}=0.
$$
or
\bb
I\equiv i\int_{\G(t)}u^*du=\s=\const. \label{inv}
\ee
The integral in this formula geometrically means the doubled
square of the domain bounded by the curve $\G(t)$. The
necessary condition of the boundedness of the first
correction term is invariance of this square. It is the
adiabatic invariant.
\par
The solution of the equation for the first correction term
can be written in the form

\bb
\stackrel{1}{U}(t_1,t)=c_1\pt_{t_1}\stackrel{0}{U}+
{-i\phi'\over\pt_E S}\pt_E\stackrel{0}{U}+\stackrel{1}{R}(t_1,t).
\label{u1-form}
\ee
\par
Where
$$
\stackrel{1}{R}(t_1,t)=u_1(t_1,t)\int_0^{t_1}
dt'\big(\pt_t\stackrel{0}{U}(t',t)u_2^*(t',t)+
\pt_t\stackrel{0}{U}\!^*(t',t)u_2(t',t)\big) \,-\,
$$
$$
\,-\, u_2(t_1,t)\int_0^{t_1}
dt'\big(\pt_t\stackrel{0}{U}(t',t)u_1^*(t',t)+
\pt_t\stackrel{0}{U}\!^*(t',t)u_1(t',t)\big).
$$
To evaluate the second correction we use the property of
even and odd of the real and imaginary parts of the first
correction term. These properties of the first and the second
terms in the formula (\ref{u1-form}) can be obtained from analysis
of the leader term of the asymptotics. The last term of this sum
contains the even function $\Re(\stackrel{1}{R})$ and
the odd function $\Im(\stackrel{1}{R})$ with respect to variable
$t_1$.
\par
Constructing of the periodical second correction of the
asymptotics leads us to determining of the function $\phi(t)$.
The necessary condition of this property for the second
correction is $$ \int_0^Tdt_1\bigg(\stackrel{2}{F}
\pt_{t_1}\stackrel{0}{U}\!^* +\stackrel{2}{F}\!^*
\pt_{t_1}\stackrel{0}{U}\bigg)=0. $$ Where  $\stackrel{2}{F}$
is the right-hand side of the equation (\ref{u2}).
\par
Let us separately evaluate the integral depending on square
of the term $\stackrel{1}{U}$

$$
J=\int_0^Tdt_1\big[-(2|\stackrel{1}{U}|^2\stackrel{0}{U}+
\stackrel{1}{U}\!^2\stackrel{0}{U}\!^*) \pt_{t_1}\stackrel{0}{U}\!^*
\,-\, (2|\stackrel{1}{U}|^2\stackrel{0}{U}\!^*+
\big(\stackrel{1}{U}\!^*\big)^2\stackrel{0}{U})
\pt_{t_1}\stackrel{0}{U} \big].
$$
Integrating by parts and substituting of the expressions
for $\pt_{t_1}\stackrel{1}{U}$ leads us to the following formula
$$
J=\int_0^Tdt_1\big[\big(-i\pt_t\stackrel{0}{U}
-i\phi'\pt_{t_1}\stackrel{0}{U}\big) \pt_{t_1}\stackrel{1}{U}\!^*
+\big(-i\pt_t\stackrel{0}{U}\!^*
+i\phi'\pt_{t_1}\stackrel{0}{U}\!^*\big) \pt_{t_1}\stackrel{1}{U}
\big]
$$
Then we obtain
$$
\pt_t\int_0^Tdt_1\big[i\stackrel{1}{U}\!^*\pt_{t_1}\stackrel{0}{U}
\,-\,i\stackrel{1}{U}\pt_{t_1}\stackrel{0}{U}\!^*\big].
$$
Substituting of $\stackrel{1}{U}$ in this formula and
using of even properties gives us
$$
\pt_t\int_0^Tdt_1\,{\phi'\over\pt_E
S}\big(\pt_E\stackrel{0}{U}\pt_{t_1}\stackrel{0}{U}\!^*\big)=0
$$
or \bb {\phi'\over \pt_E S}\pt_E I \equiv\phi_1=\const.
\label{phi'} \ee
This elegant form of the equation for the phase shift was
obtained by F.J.Bo\-ur\-land и R.Haberman \cite{Bour-Hab2} for
the solutions of the second order differential equation.
\par
The formulas (\ref{u0s}), (\ref{period}), (\ref{inv}) and
(\ref{phi'}) determine the family of the fast oscillating
asymptotic solutions of the equation (\ref{sh}). The
solutions of this family contain three parameters. The
constant $\s$ determines the trajectory on the phase plane.
The constant $\phi_1$ in the formula (\ref{phi'}) and initial
condition $\phi|_{t=\const}=\phi_0$ determine the phase
shift. The periodic with respect to $t_1$ leader term of the
asymptotics depends on three constants. Specificity  of this
fact is explained by the large values of the time variable.
So to determine the phase shift correctly we should take into
account one condition for the first correction term of the
constructed asymptotics. The period $T$ is arbitrary nonzero
constant.

\subsection{The confluent asymptotic solution}

\par
In this section we choose the values of the parameters when
the constructing fast oscillating solution degenerates at the
moment $t=t_*$. The algebraic solution constructed above
degenerates at the moment $t=t_*$ in the point $U_*$. And
internal asymptotic solution has as the leader term the
constant $U_*$. It leads to the fact that the fast oscillating
solution  degenerates at the point $U_*$ and it allows us to
match constructed solutions in the different domains. This
reason allows us to evaluate the value $E(t_*)=E_*$
$$
E_*={1\over2}|U_*|^4-t|U_*|^2-(U_* + U_*^*).
$$
By substituting these values we obtain
$E_*={3\over4}\left({1\over2}\right)^{1/3}$.
At the point $U_*$ all three roots of the equation
$$
2y^3+(2E+t_*^2)y^2 + 2t_*y+1=0
$$
coalesce. It means that the integrand in the equation
(\ref{u0s}) has the singularity of the order of $3/2$ and
elliptic integral in (\ref{u0s}) degenerates. The obtained
values of the constants $E_*$ and $t_*$ allow us to evaluate
the adiabatic invariant $\s$. To obtain this value it is
necessary to evaluate the integral
\bb
 I_*=\int_{\G(t_*)}u^* du=\s_*. \label{sigma}
\ee
As result the equation (\ref{inv}) defines the function $E(t)$
for our confluent solution. Thus to obtain the leader term of
the asymptotics (\ref{as4}) completely it is necessary to
determine the phase shift $S_0$ in the formula ({\ref{u0s}).
It can be evaluated by matching of the asymptotics
(\ref{as4}) and internal asymptotics from previous section.

\subsection{The domain of suitability of confluent asymptotic
solution as $t\to t_*-0$}

To match the internal asymptotic expansions with the fast
oscillating expansion it is necessary to know where the last
one is valid as $t\to t_*$.  From the formula (\ref{period})
it is easy to see that the function $S'$ is defined by the
value inverse to an integral
$$
K=\int_{\G(t)}{dy\over\sqrt{2y^3+(2E+t^2)y^2+2ty+1}}.
$$
So to obtain the domain of suitability of the fast oscillating
solution it is necessary to determine the order of
singularity of the leader term of the asymptotics of this
integral as $t\to t_*-0$.
\par
The integral $K$ contains two parameters $t$ and $E$. The
function $E=E(t)$ is defined by equation (\ref{inv}). In the
neighborhood of the point $t_*$ we use the following notation
$$
t=t_*+\mu,\quad E=E_*+\d.
$$
Let us change the variable in the integral
$$
y=-\bigg({1\over2}\bigg)^{1/3}+\z.
$$
Then the integral $K$ looks like
\ba
K=\int_{\G(t)}d\z\bigg(2\z^3+
\z^2\bigg[2\d+6\bigg({1\over2}\bigg)^{1/3}\mu+\mu^2\bigg]+
\\
\z\bigg[-2^{5/3}\d-4\mu-2^{2/3}\mu^2\bigg]+\\
\bigg[2^{1/3}\d+3\bigg({1\over2}\bigg)^{1/3}\mu-2^{2/3}\mu+
\bigg({1\over2}\bigg)^{2/3}\mu^2\bigg]\bigg)^{-1/2}. \ea
\par
The curve $\G(t)$ as $t=t_*$ passes through a point $\z=0$.
At this point the integrand at the moment $t=t_*$ has
nonintegrable singularity and integral diverges.
\par
The integrand and curve $\G(t)$ give us the information that
to evaluate the main term of the asymptotics as $t\to t_*-0$
 of this integral it is enough to evaluate the integral only on
the small segment of the curve $\G(t)$ near the point $\z=0$.
\par
To make this evaluations it is convenient to pass to the
integral on a small segment of the real axis and to write the
complex variable $\z$ in the form
$$
\z=\xi+i\eta,\quad  \xi,\eta\in {\hbox{R}}.
$$
In a small neighborhood of point $\z=0$
the imaginary part of the $\z$ can be expressed by real part
by formula
$$
\eta_{\pm}(\xi)=
$$
$$
\pm\sqrt{2^{1/3}+2^{2/3}\xi-\xi^2+\mu-
2^{1/3}\sqrt{1+3\bigg({1\over2}\bigg)^{1/3}\mu+2^{1/3}\d+
2^{-2/3}\mu^2+2^{4/3}\xi}}.
$$
Now if we assume that $\xi$, $\mu$ and  $\d$ are small then
we obtain \ban \eta_{\pm}^2=
-2^{-2/3}\xi^3+2^{-1/3}\bigg[3\bigg({1\over2}\bigg)^{1/3}\mu+
2^{1/3}\d\bigg]\xi+\bigg[-\mu/2-2^{1/3}\d\bigg]+
\nonumber\\
O(\xi^4)+O(\xi^2\mu)+ O(\xi^2\d)+
O(\xi\d)+O(\xi\mu)+O(\mu^2)+O(\d^2)+O(\mu\d). \label{eta-as}
\ean
The function $\d(\mu)$  stays undetermined in this
formula.
\par
The dependence $E(t)$ and hence $\d(\mu)$ can be obtained from
(\ref{inv}). Let us write this equation by integrals on
segments
\bb
\int_{x_l}^{x_+}dx\sqrt{t-x^2+\sqrt{t^2+2E+4x}}-
\int_{x_l}^{x_-}dx\sqrt{t-x^2-\sqrt{t^2+2E+4x}}=I_*. \label{E}
\ee
Where $x_l=-(t^2+2E)/4$, $x_+$ is  the right real root of the
equation
\bb x^4-2tx^2-4x-2E=0, \label{p4} \ee
and $x_-$ is the left one.
\par
To construct the asymptotics of the function $E(t)$ as $t\to
t_*-0$ it's necessary to write the asymptotics of the roots
of the equation (\ref{p4}).  This equation in the investigated
domain has two real roots and two complex conjugate roots. So
we write

$$
(x-x_+)(x-x_-)((x-m)^2+n^2)=0.
$$
Where $x_-$, $x_+$,  $m$ и $n$ are real functions  with
respect to $t$. They satisfy to Vieta equations

\ba
x_-x_+(m^2+n^2)=-2E,\\ -(x_-+x_+)(m^2+n^2)-2mx_-x_+=-4,\\
m^2+n^2+x_-x_++2m(x_-+x_+)=-2t,\\x_-+x_++2m=0.
\ea
\par
We seek a solution of these system of equations in the form
$$
x_-=U_*+\a,\quad x_+=-4U_*+\b,\quad  m=U_*+m_1.
$$
The substituting of the values $U_*$, $E_*$, $t_*$ leads us to
\ba
2^{2/3}[3m_1^2-n^2-\mu]-2m_1[-m_1^2+n^2+\mu]=0,\\
2^{1/3}[6m_1^2-2n^2-\mu]+2^{3/5}m_1[\mu-3m_1^2-2n^2]+
\\
+[m_1^2+n^2][-2\mu+3m_1^2-n^2]=-2\d. \ea
\par
The asymptotics of the solution for this system of the
equation is seeking in the form \ba
m_1=\m_1\sqrt{-\mu}+\m_2\mu+\dots,\\
n=\n_1\sqrt{-\mu}+\n_2\mu+\dots,\\
\d=\d_1\mu+\d_2(-\mu)^{3/2}+\dots. \ea By substituting into
the system we get
$$
3\m_1^2-\n_1^2=1,\quad \d_1=-\bigg({1\over2}\bigg)^{2/3}.
$$
Thus in the first corrections we did not determine a relations
between $\m_1$ and $\n_1$. This relation can be determine  from integral
equation (\ref{E}).
\par
Let us construct the main term of the asymptotics of the
integral $K$. In the small neighborhood of the point $\z=0$
the curve $\G(t)$ intersects of the real axis. So the integral
over a small neighborhood of the point $\z=0$ can be
represented in the form
 \ba K_0\sim\int_\a^{\xi_{max}}{\pt_\xi\z_+(\xi)
d\xi\over\sqrt{2\z_+^3-2[-2^{4/3}+2]\mu\z_+
+\a_+(-\mu)^{3/2}}}+\\
\int_\a^{\xi_{max}}{\pt_\xi\z_-(\xi)
d\xi\over\sqrt{2\z_-^3-2[-2^{4/3}+2]\mu\z_-
+\a_-(-\mu)^{3/2}}}, \ea where $\xi_{max}$ is the null of the
radicand, $\z_{\pm}=\xi+i\eta_{\pm}(\xi)$, $\a_{\pm}$ are
constants. The parameter $\a$ is chosen such that
$\a(-\mu)^{-1/2}\gg1$ and $\a\ll1$ as $\mu\to-0.$
\par
The asymptotics of the function $\eta_{\pm}$ and the explicit
formula for the main term of the integral $I_0$ and change of
the variable $\xi=(-\mu)^{1/2}\nu$ in integral
 lead us to the formula
$$
K_0=O((-\mu)^{-1/4}).
$$
Therefore
$$
S'=O((-\mu)^{1/4}),\quad \hbox{as}\ \mu \to -0.
$$
By the similar evaluation one can show that
$$
\pt_E I=O((-\mu)^{-1/4}),\quad  \pt_E S=O((-\mu)^{1/4}),
{\hbox{as}}\,\,\mu \to -0.
$$
Therefore  $\phi'=O((-\mu)^{1/2})$.
\par
By this way we obtain the formula
$$
\pt_t\stackrel{0}{U} =  O(\mu^{-1})\quad
{\hbox{as}}\,\,\mu\to-0.
$$
\par
Using $u_1$ and $u_2$  --- solutions  of the homogeneous
equation one can obtain that
$$
\stackrel{1}{U}=O((-\mu)^{-3/2}).
$$
\par
As a result we obtain that  the fast oscillating asymptotics
is valid as $(t_*-t)\ve^{-2/3} \gg 1$. It is easy to see that
the domains, where the constructed asymptotic solutions are
valid, are intersected.
\par

\subsection{Matching of the asymptotics}

\par
Let us match the fast oscillating asymptotics and internal
asympotics
 (\ref{infs}), (\ref{weia}) into the transition layer.
\par
The main term of the fast oscillating asymptotics is defined as
implicit function by the formula
$$
t_1=-iS'\int^{U_*+W}_{y\in\G(t)}{dy\over -|y|^2y-ty+1}.
$$
\par
This integral can be represented in the form as
$t\to t_*+0$
$$
t_1=-iS'\bigg(\int_{\G(t_*)}^{U_*+W}{dy\over
-|y|^2y-t_*y+1+O(t-t_*)y} \,+\,
$$
$$
\int_{\d\S_{\G}}^{U_*+W}
dy\wedge dy^*\pt_{y^*}\bigg[{1\over
-|y|^2y-t_*+1+O(t-t_*)y}\bigg]\bigg).
$$
Where $\d{\S_{\G}}$ is a figure bounded by the curves $\G(t)$
and $\G(t_*)$.
\par
The function $U_*+W$ is periodic. In the case the following
inequalities are
$$
{\d\s_{\G}\over W^2}\ll1\quad \hbox{(where}\,\, \d\s_{\G} \hbox{--
square of the figure }\,\,\d\S_{\G}) \quad {t-t_*\over W^{7/2}}\ll1,
$$
we obtain the following representation in  some neighborhoods
of the point   $t_{(k)}$
$$
\stackrel{0}{U}\sim U_*+W(t_{(k)}),\quad
$$
Where
$$
W(t_{(k)})={-2\over (t_{(k)}-iU_*)^2}.
$$
The argument of the function $W$ can be represented in the form
$$
{t_1\over S'}\sim
{S(t_{(k)})+S'(t_{(k)})(t-t_{(k)})+o(t-t_{(k)})\over
\ve(S'(t_{(k)})+O(S''(t_{(k)})(t-t_{(k)}))} \,+\,
{\phi(t_{(k)})\over S'(t_{(k)})},\quad {\hbox{as}}\,\,\,\,
t\to t_{(k)}.
$$
By matching with the main term of the internal asymptotics we
obtain
$$
{S(t_{(k)})\over \ve S'(t_{(k)})} \,+\, {t-t_{(k)}\over \ve}
\,+\, {\phi(t_{(k)})\over S'(t_{(k)})}\sim \th_k.
$$
\par
This matching with the second internal asymptotics gives us
that this fast oscillating solution also can be matched with
the intermediate expansion. The constants $\phi_0$ and
$\phi_1$ evaluate by this way. These constants define the
phase shift of the fast oscillating solution. But this kind of
evaluations are not made in this work.
\par
We are grateful to L.A. Kalyakin, B.I. Suleimanov for
stimulating discussions. And we also thank  N. Enikeev for
his helpful and comments.

\end{document}